\pdfoutput=1
\documentclass[11pt]{article}
\usepackage[a4paper,margin=27mm]{geometry}
\usepackage[T1]{fontenc}
\usepackage{lmodern,microtype,amsmath,amssymb,amsthm,mathtools}
\usepackage{enumitem,xcolor}
\usepackage[numbers,square,sort&compress]{natbib}
\usepackage[colorlinks=true,linkcolor=blue!45!black,citecolor=blue!45!black,urlcolor=blue!45!black]{hyperref}
\usepackage[capitalise]{cleveref}
\hypersetup{pdftitle={Kolmogorov Obstacle Problems: Euclidean Free Boundary Regularity under a Cone Condition},pdfauthor={Junkee Jeon, Takwon Kim, Jehan Oh}}
\newtheorem{theorem}{Theorem}[section]
\newtheorem{proposition}[theorem]{Proposition}
\newtheorem{lemma}[theorem]{Lemma}
\newtheorem{corollary}[theorem]{Corollary}
\theoremstyle{definition}

\theoremstyle{remark}
\newtheorem{remark}[theorem]{Remark}
\numberwithin{equation}{section}
\newcommand{\R}{\mathbb R}

\title{\bfseries Kolmogorov Obstacle Problems:\\ Euclidean Free Boundary Regularity\\ under a Cone Condition}
\author{Junkee Jeon\thanks{Department of Applied Mathematics, Kyung Hee University, Yongin 17104, Republic of Korea. E-mail: junkeejeon@khu.ac.kr}
\and Takwon Kim\thanks{School of Mathematics, Statistics and Data Science, Sungshin Women's University, Seoul 02844, Republic of Korea. E-mail: takwon@sungshin.ac.kr}
\and Jehan Oh\thanks{Corresponding author. Department of Mathematics, Kyungpook National University, Daegu 41566, Republic of Korea. E-mail: jehan.oh@knu.ac.kr}}
\date{}

\begin{document}
\maketitle

\begin{abstract}
We study obstacle problems for Kolmogorov-type operators that diffuse in one variable and transport in another. Such problems describe American options on running averages. Under a cone estimate, a strong monotonicity condition on the solution in the diffusive direction, we show that the free boundary is a graph in the diffusive direction with Hölder continuous first derivatives in the Euclidean sense, and that every free boundary point is regular. The key observation is that a Euclidean Lipschitz free boundary is flat at the intrinsic scale. For the American fixed-strike Asian option we verify the cone estimate, so that its exercise boundary has the same regularity.
\end{abstract}



\section{Introduction and main results}\label{sec:intro}
Kolmogorov operators combine diffusion in some variables with transport in the others. In this paper we study obstacle problems for operators of the form
\begin{equation}\label{eq:K}
 \mathcal{K}u=a\,u_{xx}+\beta\,u_x+x\,u_y+u_t-ru,
\end{equation}
which diffuse in $x$ and transport in $y$ with velocity $x$. For $a=\frac12$ and $\beta=r=0$, the operator $\mathcal{K}$ is the backward operator of a Brownian motion together with its time integral, whose transition density was computed by Kolmogorov~\cite{Kol34}. Operators of this type are the prototypes of kinetic Fokker--Planck operators. In mathematical finance they arise whenever the state of a model contains the running integral of a diffusion, as for options written on the average of a stock price~\cite{BPV01}. Optimal stopping for such models, for instance the pricing of American Asian options, leads to obstacle problems for $\mathcal{K}$. So does singular control with a degenerate diffusion, as in the consumption-adjustment problem of~\cite{HJO}, where the risk premium is mean reverting and the two state variables are driven by the same Brownian motion.

Although $\mathcal{K}$ diffuses only in $x$, it is hypoelliptic. Indeed, the commutator of $\partial_x$ with the transport field $Y=x\partial_y+\partial_t$ is $\partial_y$, so that $\partial_x$, $Y$ and $\partial_y$ span all directions~\cite{Hor67}. The natural geometry of $\mathcal{K}$ is not Euclidean. For constant $a$, the operator $a\partial_{xx}+Y$ is invariant under the translations and homogeneous of degree two under the dilations
\[
 (x,y,t)\mapsto(x_P+x,\ y_P+y+x_Pt,\ t_P+t)\qquad\text{and}\qquad(x,y,t)\mapsto(\rho x,\ \rho^3y,\ \rho^2t),
\]
respectively~\cite{LP94}. The translations follow the transport. Accordingly, an intrinsic cylinder of size $\rho$ has extent $\rho$ in the diffusive variable and $\rho^2$ in time, and its extent in the transport variable is $\rho^3$ about the line along which its center is transported, see \eqref{eq:cylinder} below. After the obstacle is subtracted, the obstacle problem takes the normalized form
\begin{equation}\label{eq:obstacle}
 u\ge0\qquad\text{and}\qquad\mathcal{K}u=f\chi_{\{u>0\}},
\end{equation}
where $u$ is the difference between the value function and the obstacle, and the positive function $f$ is determined by the obstacle. Solutions of \eqref{eq:obstacle} have bounded second derivatives in the intrinsic sense, namely $u_{xx},Yu\in L^\infty_{\mathrm{loc}}$, and they grow quadratically at free boundary points with respect to the intrinsic scale~\cite{FNPP10}. Where the free boundary $\Gamma=\partial\{u>0\}$ is a surface in $\R^3$, its position relative to the diffusion matters. Near a non-characteristic point, $\Gamma$ is a graph in the diffusive direction and the diffusion acts across it. Where $\Gamma$ is not such a graph, for instance where the diffusion is tangential to $\Gamma$, the behavior of $u$ near $\Gamma$ is governed by the transport. Where the transport is tangential as well, as on the grazing set of kinetic boundary problems with in-flow conditions, solutions lose regularity~\cite{KW26b}.

For the Laplace and heat operators, the regularity theory of free boundaries in obstacle problems is classical. Caffarelli~\cite{Caf77,Caf98} divided free boundary points into regular and singular ones according to their blow-ups, which are half-space solutions at regular points and quadratic polynomials at singular points. Near regular points the free boundary is continuously differentiable~\cite{Caf77}, and by the work of Kinderlehrer and Nirenberg~\cite{KN77} it is then analytic for the model problem. The parabolic theory was developed in~\cite{Caf77,CPS04}, and we refer to the monograph~\cite{PSU12} for both settings. In optimal stopping the payoff often provides additional structure. For American options on several assets with a convex payoff, Laurence and Salsa~\cite{LS09} used qualitative properties of the exercise region to prove that its boundary is $C^\infty$.

For hypoelliptic operators the theory is much less developed. On Carnot groups of step two, Danielli, Garofalo and Petrosyan~\cite{DGP07} established the $C^{1,\alpha}$ regularity of regular free boundaries for the sub-Laplacian. For Kolmogorov operators, the obstacle problem was introduced in connection with American Asian options~\cite{DFPP08,Pas08,MP09}. Frentz, Nyström, Pascucci and Polidoro~\cite{FNPP10} proved the optimal regularity of solutions in the intrinsic sense, see also~\cite{NPP10} for the behavior near the initial state. The free boundary has been studied only recently. Bowman~\cite{Bow26} showed that $u_y$ is locally bounded, which is more than the intrinsic scaling gives. With this estimate he derived an almost monotonicity formula of the type introduced by Weiss~\cite{Wei99} and showed that blow-ups solve the parabolic obstacle problem, so that free boundary points are again regular or singular as in~\cite{CPS04}. Under a thickness condition on the contact set, he proved that the free boundary is a graph $x=\Phi(y,t)$ in the diffusive direction, where $\Phi$ is Hölder continuous of order $\frac12$, and that it is differentiable in the intrinsic sense. Nyström~\cite{Nys26b} then showed that, near every regular point, the free boundary is a non-characteristic surface of class $C^{1,\alpha}$ in the intrinsic sense. His proof uses boundary Harnack inequalities in asymptotically cylindrical domains~\cite{Nys26a}, which build on the boundary estimates of~\cite{CNP12,CNP13,NP16,LN22}.

Two questions are left open by these results, and both matter in applications. The first is which free boundary points are regular. The thickness condition rules out singular points, but it is a condition on the unknown contact set, and it is difficult to verify for a given problem. The second is the regularity of the free boundary in the Euclidean sense. Intrinsic regularity describes the free boundary in intrinsic cylinders, in which $x$ varies by $\rho$ while $y$ and $t$ vary only by $O(\rho^2)$. In applications, however, the free boundary is an exercise boundary, that is, a function of calendar time and of the state, and its Euclidean regularity is the natural question.

In this paper we show that the monotonicity of the data in the application answers both questions. In the American Asian option the value function is nondecreasing in the stock price, which is the diffusive variable, while the obstacle depends only on the average. For data that are monotone in the diffusive direction in this way, and under a cone estimate for the solution, the free boundary is a $C^{1,\gamma}$ graph in the Euclidean sense, and every free boundary point is regular and non-characteristic, see \cref{thm:A}. The cone estimate is a strong assumption, and its verification is the main point of the application. For the Asian option we obtain it from the structure of the payoff, so that its exercise boundary has the same regularity, see \cref{cor:asian}.

We now state the setting and the results. Let $Q\subset\R^3$ be an open set with coordinates $(x,y,t)$. We consider the operator \eqref{eq:K} with coefficients $a,\beta\in C^\infty(\overline Q)$ and a constant $r\ge0$, and we assume
\begin{equation}\label{eq:H}
 a\ge a_0>0\qquad\text{in }\overline Q.
\end{equation}
A transport velocity $c(x)$ with $c'\neq0$ reduces to this case through the change of variables $x\mapsto c(x)$. A solution of \eqref{eq:obstacle} in $Q$ is a function $u\in C(Q)$ with $u_x,u_{xx},Yu\in L^\infty_{\mathrm{loc}}(Q)$ such that \eqref{eq:obstacle} holds almost everywhere in $Q$, where $f\in C^\infty(\overline Q)$ satisfies $f\ge f_0>0$. Here the derivatives are taken in the sense of distributions. This regularity is optimal in the intrinsic sense~\cite{FNPP10}. By Hörmander's theorem~\cite{Hor67}, $u$ is smooth in the open set $\{u>0\}$. The free boundary is $\Gamma=\partial\{u>0\}\cap Q$.

The intrinsic cylinders follow the transport. For $P=(x_P,y_P,t_P)\in\R^3$ and $\rho>0$ we set
\begin{equation}\label{eq:cylinder}
 Q_\rho(P)=\bigl\{(x,y,t):\ |x-x_P|<\rho,\ |t-t_P|<\rho^2,\ |y-y_P-x_P(t-t_P)|<\rho^3\bigr\}.
\end{equation}
With the maps
\begin{equation}\label{eq:Lambda}
 \Lambda_{P,\rho}(x,y,t)=\bigl(x_P+\rho x,\ y_P+x_P\rho^2t+\rho^3y,\ t_P+\rho^2t\bigr)
\end{equation}
we have $\Lambda_{P,\rho}(Q_R(0))=Q_{R\rho}(P)$ for $R>0$, and a direct computation gives, for every smooth function $w$,
\begin{equation}\label{eq:Lambdaderiv}
 \partial_{xx}(w\circ\Lambda_{P,\rho})=\rho^2(\partial_{xx}w)\circ\Lambda_{P,\rho}\qquad\text{and}\qquad Y(w\circ\Lambda_{P,\rho})=\rho^2(Yw)\circ\Lambda_{P,\rho}.
\end{equation}
A point $P\in\Gamma$ is regular if every locally uniform limit of the rescalings $u_\rho=\rho^{-2}u\circ\Lambda_{P,\rho}$ along a sequence $\rho\searrow0$ is a half-space solution $\frac{f(P)}{2a(P)}\bigl((\epsilon x)^+\bigr)^2$ with $\epsilon\in\{-1,1\}$. We call $P\in\Gamma$ non-characteristic if, near $P$, the set $\Gamma$ is the graph $\{x=\Phi(y,t)\}$ of a Lipschitz function. Then the diffusive direction is transversal to $\Gamma$ at $P$. This is stronger than the intrinsic notion of~\cite{Nys26a,Nys26b}, where $\Phi$ is only required to be Lipschitz with respect to the intrinsic distance.

We call a solution of \eqref{eq:obstacle} diffusive-monotone if $u_x\ge0$ almost everywhere in $Q$. Our main result states that a cone estimate, which implies diffusive monotonicity, makes the free boundary a graph of class $C^{1,\gamma}$ and every free boundary point regular.
\begin{theorem}\label{thm:A}
Let $u$ be a solution of \eqref{eq:obstacle} under \eqref{eq:H} that satisfies the cone estimate
\begin{equation}\label{eq:cone}
 |u_y|+|u_t|\le C_0\,u_x\qquad\text{in }Q\cap\{u>0\}
\end{equation}
with a constant $C_0>0$. Then the following statements hold for every open set $Q'\Subset Q$.
\begin{enumerate}[label=\textup{(\roman*)}]
\item Near every point of $\Gamma$, the set $\Gamma$ is the graph $\{x=\Phi(y,t)\}$ of a function $\Phi$ with
\[
 |\Phi(y,t)-\Phi(y',t')|\le C_0\max\bigl\{|y-y'|,|t-t'|\bigr\}.
\]
\item Every point of $\Gamma\cap Q'$ is regular and non-characteristic. Moreover, there is $\varepsilon>0$ such that $u$ is comparable to the square of the distance to $\Gamma$ at the points of $\{u>0\}\cap Q'$ whose distance to $\Gamma$ is less than $\varepsilon$.
\item There is $\gamma\in(0,1)$ such that, near every point of $\Gamma\cap Q'$, the function $\Phi$ in~(i) is of class $C^{1,\gamma}$.
\end{enumerate}
\end{theorem}

The cone estimate \eqref{eq:cone} is the quantitative form of diffusive monotonicity. It says that $u$ is nondecreasing along every direction $(1,\mu_1,\mu_2)$ with $|\mu_1|,|\mu_2|\le C_0^{-1}$, and statement~(i) follows from this by the classical argument for cones of monotonicity~\cite{Caf77,PSU12}. In particular $u_x\ge0$ in $\{u>0\}$, and $u_x=0$ almost everywhere in the contact set, because $u$ attains its minimum there and $x\mapsto u(x,y,t)$ is continuously differentiable for almost every $(y,t)$. Hence every solution that satisfies \eqref{eq:cone} is diffusive-monotone. Diffusive monotonicity alone makes the contact set the region below a graph in the diffusive direction, but it does not control the regularity of this graph. In the application both conditions come from the structure of the payoff, and the cone estimate requires a proof.

The cone estimate is a strong assumption on the solution. By statement~(i) and \cref{lem:flat}, in every small intrinsic cylinder $Q_\rho(P)$ centered at a free boundary point $P$, each section of the contact set in the diffusive direction is an interval of length at least $\rho-C\rho^2$. Hence the contact set is thick in the sense of~\cite{Bow26,Nys26b} at every free boundary point. Thus \cref{thm:A} does not weaken the thickness condition, and for the model problem statement~(ii) is contained in the results of~\cite{Bow26}. Our proof of~(ii) is different and short, and it applies to variable coefficients and a variable function $f$. The point of the cone estimate is that it is a condition on derivatives of the solution that can be verified from the structure of a given problem, as we do for the Asian option, while the thickness condition concerns the unknown contact set. It is a hypothesis on the solution and not on the data. The monotonicity of the data in the application is what makes it plausible, and its proof for the Asian option occupies most of \cref{sec:asian}. Under the cone estimate there are no singular points, at which blow-ups would be quadratic polynomials~\cite{Bow26}, and the comparability of $u$ with the squared distance to $\Gamma$ combines the optimal growth with nondegeneracy at every point.

Statement~(iii) gives the existence and the Euclidean Hölder continuity of the derivatives $\Phi_y$ and $\Phi_t$. This is the main difference from the intrinsic differentiability of~\cite{Bow26} and the intrinsic $C^{1,\alpha}$ regularity of~\cite{Nys26b}. For one diffusive variable the intrinsic normal of~\cite{Nys26b} takes only the values $\pm1$, and the intrinsic $C^{1,\alpha}$ regularity of $\Gamma$ near $P$ established in~\cite[Theorem~1.1]{Nys26b} amounts to the flatness bound
\begin{equation}\label{eq:intrinsic}
 |\Phi(y',t')-\Phi(y,t)|\le C\Bigl(\bigl|y'-y-(t'-t)\Phi(y,t)\bigr|^{\frac13}+|t'-t|^{\frac12}\Bigr)^{1+\alpha}
\end{equation}
for points of $\Gamma$ near $P$. This bound allows $\Phi$ to be merely Hölder continuous of order $\frac{1+\alpha}{2}$ in $t$ and of order $\frac{1+\alpha}{3}$ along the transport, and it contains no information about derivatives of $\Phi$. For example, the Lipschitz graph $x=|y|$ satisfies \eqref{eq:intrinsic} with $\alpha=1$ at every point by \cref{lem:flat}, but it is not differentiable at $y=0$. Whether free boundaries that are not differentiable in the Euclidean sense occur for Kolmogorov obstacle problems is not known. Under the cone estimate they do not occur. Conversely, a Euclidean Lipschitz graph in the diffusive direction is flat of order two at the intrinsic scale by \cref{lem:flat}, so that statement~(i) alone yields \eqref{eq:intrinsic} with exponent one in place of $\alpha$. Statement~(iii) therefore answers a question that intrinsic regularity leaves open. \Cref{thm:A} allows variable coefficients and a variable function $f$, whereas the results of~\cite{Bow26,Nys26b} concern the model problem with $f=1$.

\Cref{thm:A} concerns one diffusive variable, whereas the results of~\cite{Bow26,Nys26b} hold in every dimension. For operators with several diffusive variables $x\in\R^m$ and transport $x\cdot\nabla_y$, the analogous cone estimate makes $\Gamma$ a Lipschitz graph $x_m=\Phi(x',y,t)$, and two steps of our proof then fail. First, the variable $x'$ ranges over an interval of length $\rho$ in an intrinsic cylinder of size $\rho$, so a Lipschitz graph deviates from the hyperplane $\{x_m=x_{P,m}\}$ by $O(\rho)$ and not by $O(\rho^2)$ as in \cref{lem:flat}. The blow-ups then retain their dependence on $x'$. Since the cone estimate removes the dependence on $y$ and $t$ in the limit, they are global solutions of the elliptic obstacle problem in $x$ whose contact sets contain a cone, and \cref{lem:liouville} would have to be replaced by Weiss's monotonicity formula and the classification of homogeneous global solutions~\cite{Caf98,PSU12}, which are available for these elliptic blow-ups. We therefore expect that statement~(ii) extends to several diffusive variables, although we have not carried this out. Second, in the proof of statement~(iii) the drift $\beta u_x$ is removed by an exponential conjugation. With several diffusive variables this requires that $a^{-1}\beta$ be a gradient field, up to the divergence of $a$, and otherwise the comparison principle of~\cite{LN22}, which is stated for operators without drift, would have to be extended. This is the actual obstruction. Since the application in \cref{sec:asian} has one diffusive variable, we do not pursue the general case. In summary, \cref{thm:A} allows variable coefficients and a variable function $f$, but only one diffusive variable, and it replaces the thickness condition of~\cite{Bow26,Nys26b} by the cone estimate. In return it gives Euclidean instead of intrinsic regularity, and it does so at every free boundary point, whereas the thickness condition at a point controls the free boundary only near this point.

For the American Asian option the hypotheses of \cref{thm:A} can be verified away from the initial and terminal times. Let $S$ be a geometric Brownian motion with volatility $\sigma>0$, interest rate $r\ge0$ and dividend yield $q\ge0$ under the pricing measure, let $I_t=\int_0^tS_\tau\,d\tau$ and $A_t=I_t/t$, and let $K>0$ and $T>0$ be the strike and the maturity. The value of the American option with payoff $(A_t-K)^+$ is
\[
 V(t,S,I)=\sup_\theta\mathbb{E}\Bigl[e^{-r(\theta-t)}\Bigl(\frac{I_\theta}{\theta}-K\Bigr)^+\Bigr],
\]
where the supremum is taken over stopping times $\theta$ with values in $[t,T]$, and the processes start from $S_t=S$ and $I_t=I$.
\begin{corollary}\label{cor:asian}
There is a positive function $S_*$ on $\{0<t<T,\ A>K\}$ such that the exercise region in this set is $\{S\le S_*(t,A)\}$. The function $S_*$ is continuously differentiable, and its first derivatives are locally Hölder continuous.
\end{corollary}

In \cref{cor:asian}, $S_*(t,A)$ is the critical stock price below which immediate exercise is optimal when the average is $A$. No step of the argument requires a positive interest rate, so that the case $r=0$ is included. The obstacle problem for this option was studied in~\cite{Pas08,MP09,FNPP10}, and qualitative properties of its exercise region were obtained in~\cite{DK06}. Diffusive monotonicity holds because the value function is nondecreasing in the stock price and the obstacle $A-K$ does not depend on it. The cone estimate is the main point. We obtain it from two comparisons of payoffs, which raise the average together with the strike or move the starting time forward, and from Euler's identity for the value function, which is homogeneous of degree one in the stock price, the running integral and the strike, see \cref{lem:euler,lem:asian-shift}. These comparisons bound $-u_I$ and $u_t$ from above by multiples of $Su_S$. A lower bound for $u_t$ cannot be obtained in this way, because moving the starting time forward shortens the horizon. We prove it in \cref{prop:asian-cone} with the local comparison principle of Litsgård and Nyström~\cite{LN22}, which we apply at a single scale and restate in \cref{sec:euclid}.

Let us briefly describe the ideas of the proofs. The main difficulty for Kolmogorov operators is the anisotropy of the intrinsic geometry. Blow-ups are taken with the dilations $(\rho x,\rho^3y,\rho^2t)$, so the blow-up analysis sees the free boundary only at the intrinsic scale. For graphs in the diffusive direction, however, the two geometries fit together. By \cref{lem:flat}, a Euclidean Lipschitz graph $x=\Phi(y,t)$ lies within $O(\rho^2)$ of the hyperplane $\{x=x_P\}$ in $Q_\rho(P)$. The intrinsic rescaling divides $x$ by $\rho$, so the rescaled free boundaries converge to this hyperplane, and every blow-up vanishes on a half-space. The Liouville theorem in \cref{lem:liouville} shows that such a blow-up is a half-space solution. Its proof uses only Hörmander's theorem and the invariance \eqref{eq:Lambdaderiv} of the operator, and we need neither a monotonicity formula nor a classification of all global solutions. This gives statement~(ii) of \cref{thm:A}, which therefore uses the cone estimate only through the Lipschitz graph property in statement~(i).

For statement~(iii) we show that the quotients $u_y/u_x$ and $u_t/u_x$ converge uniformly at $\Gamma$ to Hölder continuous boundary values. This suffices, because by \cref{lem:levelsets} the gradient of $\Phi$ is given by the boundary values of $-(u_y,u_t)/u_x$. The boundary values come from a decay estimate for the oscillation of the quotients at $\Gamma$, which is a boundary Harnack inequality for derivatives of $u$. These derivatives solve the equation only up to bounded right-hand sides, because $\partial_x$ does not commute with the transport and the data vary. After a multiple of $u$ is subtracted, the right-hand sides are of the size of the distance to $\Gamma$. Boundary Harnack inequalities with right-hand side are known for elliptic equations~\cite{RTL21}. The boundary Harnack inequality of~\cite{Nys26a} for Kolmogorov operators rests on a balance between forward and backward reference values, and this balance is its main difficulty. For our quotients the balance is automatic, because the denominator is comparable to the distance to $\Gamma$ and the numerators are bounded by multiples of it, see \cref{lem:linear}. Therefore we only need the local comparison principle of Litsgård and Nyström~\cite{LN22}, which we restate in \cref{sec:euclid}, and a barrier absorbs the right-hand sides, see \cref{prop:C1gamma}.

We do not know whether $\Phi$ in \cref{thm:A} is $C^\infty$ when the data are smooth. For the Laplace operator this step follows from the hodograph transform~\cite{KN77}, and for the Laplace and heat operators also from higher-order boundary Harnack inequalities~\cite{DS15,BG16,Kuk22}. For Kolmogorov operators both approaches require a Schauder estimate up to a boundary that is a graph in the diffusive variable, which does not seem to be available, see \cref{rem:higher}. At boundaries in the transport variable, such as the boundaries of spatial domains in kinetic theory, the situation is different. Kim and Weidner proved that solutions of kinetic Fokker--Planck equations with in-flow or diffuse reflection boundary conditions are Hölder continuous of order $\frac12$ up to the grazing set, which is optimal~\cite{KW26b}. They also showed that, unlike for elliptic and parabolic equations, the quotient of two solutions is in general not smooth there~\cite{KW26a}.

The paper is organized as follows. \Cref{sec:flat} proves statements~(i) and~(ii) of \cref{thm:A}. \Cref{sec:euclid} restates the boundary comparison principle of~\cite{LN22} in the form used there and in \cref{sec:asian}, proves statement~(iii) and discusses higher regularity. \Cref{sec:asian} treats the American Asian option and proves \cref{cor:asian}.

\section{Intrinsic flatness and regular points}\label{sec:flat}
We first show that the cone estimate makes the free boundary a Lipschitz graph in the diffusive direction.
\begin{proof}[Proof of \cref{thm:A}\,(i)]
Let $\mu=(\mu_1,\mu_2)$ with $|\mu_1|,|\mu_2|\le C_0^{-1}$ and put $\nu_\mu=(1,\mu_1,\mu_2)$. By \eqref{eq:cone},
\[
 \partial_{\nu_\mu}u=u_x+\mu_1u_y+\mu_2u_t\ge u_x-C_0^{-1}\bigl(|u_y|+|u_t|\bigr)\ge0\qquad\text{in }\{u>0\}.
\]
We claim that $u(z+s\nu_\mu)\ge u(z)$ whenever $s\ge0$ and the segment from $z$ to $z+s\nu_\mu$ lies in $Q$. This is clear if $u(z)=0$. If $u(z)>0$, let $s_1$ be the supremum of all $s'\in[0,s]$ such that $u>0$ on the segment from $z$ to $z+s'\nu_\mu$. On this segment $u$ is smooth and nondecreasing in the direction $\nu_\mu$. Hence $u(z+s'\nu_\mu)\ge u(z)$ for $s'<s_1$, and by continuity also for $s'=s_1$. Since $\{u>0\}$ is open, it follows that $s_1=s$, which proves the claim.

Let $P\in\Gamma$ and choose $\varepsilon>0$ such that $Q$ contains the closure of $(x_P-3\varepsilon,x_P+3\varepsilon)\times D$, where
\[
 D=\Bigl\{(y,t):\ |y-y_P|<\frac{\varepsilon}{2C_0},\ |t-t_P|<\frac{\varepsilon}{2C_0}\Bigr\}.
\]
Let $(y,t)\in D$ and put $m=\max\{|y-y_P|,|t-t_P|\}$, so that $C_0m<\varepsilon/2$. If $m>0$, the claim with $z=(x_P-C_0m,y,t)$, $s=C_0m$ and $\mu=(y_P-y,t_P-t)/s$ gives $u(z)\le u(P)=0$. Since $u(P)=0$ and $u$ is nondecreasing in $x$, in either case
\[
 u(x,y,t)=0\qquad\text{for }x_P-3\varepsilon<x\le x_P-C_0m.
\]
Now let $x_P+C_0m<x<x_P+3\varepsilon$, and take points $z_k\to P$ with $u(z_k)>0$. For large $k$ the point $(x,y,t)$ lies on a segment that starts at $z_k$ and has a direction $\nu_\mu$ as above, so $u(x,y,t)\ge u(z_k)>0$. Hence there is a number $\Phi(y,t)$ with $|\Phi(y,t)-x_P|\le C_0m$ such that, for $x_P-3\varepsilon<x<x_P+3\varepsilon$, we have $u(x,y,t)=0$ if $x\le\Phi(y,t)$ and $u(x,y,t)>0$ if $x>\Phi(y,t)$. Now let $(y,t),(y',t')\in D$ with $m'=\max\{|y-y'|,|t-t'|\}>0$. The claim with $z=(\Phi(y',t')-C_0m',y,t)$, $s=C_0m'$ and $\mu=(y'-y,t'-t)/s$ gives $u(z)\le u(\Phi(y',t'),y',t')=0$. Here the segment lies in $(x_P-3\varepsilon,x_P+3\varepsilon)\times D$, because $|\Phi(y',t')-x_P|<\varepsilon/2$ and $C_0m'<\varepsilon$. Hence
\[
 \Phi(y,t)\ge\Phi(y',t')-C_0\max\bigl\{|y-y'|,|t-t'|\bigr\}.
\]
Exchanging the two points shows that $\Phi$ satisfies the Lipschitz bound in~(i). Therefore $\Gamma\cap\bigl((x_P-\varepsilon,x_P+\varepsilon)\times D\bigr)$ is the graph of $\Phi$.
\end{proof}

The following lemma is the geometric core of the paper. A Euclidean Lipschitz graph in the diffusive direction deviates from the hyperplane $\{x=x_P\}$ only by $O(\rho^2)$ inside $Q_\rho(P)$.
\begin{lemma}\label{lem:flat}
Let $P\in\R^3$, $\rho>0$ and $L\ge0$. Let $\Phi$ be a function on the projection of $Q_\rho(P)$ to the $(y,t)$-plane such that $\Phi(y_P,t_P)=x_P$ and
\[
 |\Phi(y,t)-\Phi(y',t')|\le L\bigl(|y-y'|+|t-t'|\bigr).
\]
Then $|\Phi(y,t)-x_P|\le L(1+|x_P|+\rho)\rho^2$ for every $(x,y,t)\in Q_\rho(P)$.
\end{lemma}
\begin{proof}
Let $(x,y,t)\in Q_\rho(P)$. By \eqref{eq:cylinder},
\[
 |t-t_P|<\rho^2\qquad\text{and}\qquad|y-y_P|\le|x_P|\,|t-t_P|+\rho^3<|x_P|\,\rho^2+\rho^3.
\]
Hence $|\Phi(y,t)-x_P|\le L(|y-y_P|+|t-t_P|)\le L(1+|x_P|+\rho)\rho^2$.
\end{proof}

Since the diffusive variable is rescaled by $\rho$, \cref{lem:flat} shows that the rescaled free boundaries converge to the hyperplane $\{x=0\}$. The blow-up limits are then identified by the following Liouville theorem.
\begin{lemma}\label{lem:liouville}
Let $\bar a>0$ and $\bar f$ be constants, and let $v\in C(\R^3)$ satisfy $v_{xx},Yv\in L^\infty_{\mathrm{loc}}(\R^3)$ and
\[
 v=0\quad\text{in }\{x\le0\},\qquad\bar a\,v_{xx}+Yv=\bar f\quad\text{almost everywhere in }\{x>0\}.
\]
If $\sup_{Q_R(0)}|v|\le C_vR^2$ for all $R\ge1$, then $v=\frac{\bar f}{2\bar a}(x^+)^2$.
\end{lemma}
\begin{proof}
Let $w=v-\frac{\bar f}{2\bar a}(x^+)^2$. The function $(x^+)^2$ has the bounded second derivative $2\chi_{\{x>0\}}$ and does not depend on $y$ and $t$. Hence $w_{xx},Yw\in L^\infty_{\mathrm{loc}}(\R^3)$, and $\bar a\,w_{xx}+Yw=0$ almost everywhere in $\{x>0\}$ and in $\{x<0\}$, where $w$ vanishes. Since $\bar a\,w_{xx}+Yw$ is a locally bounded function that vanishes almost everywhere, $w$ solves $\bar a\,w_{xx}+Yw=0$ in the sense of distributions on $\R^3$. By Hörmander's theorem~\cite{Hor67}, $w$ is smooth. The same theorem gives a constant $C_*$ such that $|h_{xxx}(0)|\le C_*\sup_{Q_1(0)}|h|$ for every $h\in L^\infty(Q_1(0))$ that solves $\bar a\,h_{xx}+Yh=0$ in the sense of distributions in $Q_1(0)$. Indeed, the space of these functions is closed in $L^\infty(Q_1(0))$, because convergence in $L^\infty$ implies convergence in the sense of distributions, its elements are smooth, and the linear map $h\mapsto h|_{Q_{1/2}(0)}$ into $C^3(\overline{Q_{1/2}(0)})$ has a closed graph, since a sequence that converges in $L^\infty$ and in $C^3$ has the same limit in both spaces. By the closed graph theorem this map is bounded.

Let $P_0=(x_0,y_0,t_0)\in\R^3$ and $R\ge1$. By \eqref{eq:Lambdaderiv}, the function $h=w\circ\Lambda_{P_0,R}$ solves $\bar a\,h_{xx}+Yh=0$, and $h_{xxx}(0)=R^3w_{xxx}(P_0)$. Moreover
\[
 \Lambda_{P_0,R}(Q_1(0))=Q_R(P_0)\subset Q_{R'}(0)\qquad\text{with}\qquad R'=|x_0|+|t_0|^{\frac12}+|y_0|^{\frac13}+2R.
\]
Therefore
\[
 |w_{xxx}(P_0)|\le C_*R^{-3}\sup_{Q_{R'}(0)}|w|\le C_*\Bigl(C_v+\frac{|\bar f|}{2\bar a}\Bigr)\frac{R'^2}{R^3},
\]
and the right-hand side tends to zero as $R\to\infty$. Hence $w_{xxx}=0$, so $x\mapsto w(x,y,t)$ is a polynomial of degree at most two for every $(y,t)$. It vanishes for $x<0$, and therefore $w=0$.
\end{proof}

Together with \cref{lem:flat}, the Liouville theorem shows that all blow-ups are half-space solutions.
\begin{proof}[Proof of \cref{thm:A}\,(ii)]
Let $Q'\Subset Q''\Subset Q'''\Subset Q$ be open sets and put
\[
 M=\|u_{xx}\|_{L^\infty(Q''')}+\|Yu\|_{L^\infty(Q''')}.
\]
By~(i), its proof and the compactness of $\Gamma\cap\overline{Q''}$, there are $\rho_0\in(0,1]$ and $C_1>0$ with the following property. For every $P\in\Gamma\cap\overline{Q''}$ we have $Q_{\rho_0}(P)\subset Q'''$, and there is a function $\Phi$ as in~(i) with $\Phi(y_P,t_P)=x_P$ such that $u=0$ in $Q_{\rho_0}(P)\cap\{x\le\Phi(y,t)\}$ and $u>0$ in $Q_{\rho_0}(P)\cap\{x>\Phi(y,t)\}$. By \cref{lem:flat},
\begin{equation}\label{eq:flatP}
 |\Phi(y,t)-x_P|\le C_1\rho^2\qquad\text{for }(x,y,t)\in Q_\rho(P)\text{ and }0<\rho\le\rho_0.
\end{equation}
Decreasing $\rho_0$, we may assume that $C_1\rho_0\le\frac12$, so that $\Phi(y,t)$ lies in the $x$-slice of $Q_{\rho_0}(P)$ for every $(x,y,t)\in Q_{\rho_0}(P)$. For almost every $(y,t)$, the function $x\mapsto u(x,y,t)$ has a bounded second derivative and attains its minimum $0$ at $x=\Phi(y,t)$. Hence
\begin{equation}\label{eq:quad}
 u\le\frac M2\bigl((x-\Phi(y,t))^+\bigr)^2\qquad\text{and}\qquad|u_x|\le M\bigl(x-\Phi(y,t)\bigr)^+\qquad\text{in }Q_{\rho_0}(P)
\end{equation}
almost everywhere, and the first bound holds everywhere by continuity.

Let $P_k\in\Gamma\cap\overline{Q''}$ converge to $P$, let $\rho_k\searrow0$, and put $u_k=\rho_k^{-2}u\circ\Lambda_k$ with $\Lambda_k=\Lambda_{P_k,\rho_k}$. By \eqref{eq:Lambdaderiv}, $\partial_{xx}u_k=u_{xx}\circ\Lambda_k$ and $Yu_k=(Yu)\circ\Lambda_k$, and moreover $\partial_xu_k=\rho_k^{-1}u_x\circ\Lambda_k$. Hence
\begin{equation}\label{eq:rescaled}
 a_k\,\partial_{xx}u_k+\rho_k\beta_k\,\partial_xu_k+Yu_k-r\rho_k^2u_k=f_k\chi_{\{u_k>0\}}
\end{equation}
almost everywhere, where $a_k=a\circ\Lambda_k$, $\beta_k=\beta\circ\Lambda_k$ and $f_k=f\circ\Lambda_k$. Fix $R\ge1$. Since $\Lambda_k(Q_R(0))=Q_{R\rho_k}(P_k)$, the bounds \eqref{eq:flatP} and \eqref{eq:quad} show that, for large $k$, $u_k>0$ in $Q_R(0)\cap\{x>C_1R^2\rho_k\}$ and
\[
 0\le u_k\le\frac M2\bigl((x+C_1R^2\rho_k)^+\bigr)^2,\qquad|\partial_xu_k|\le2MR,\qquad|\partial_{xx}u_k|+|Yu_k|\le M\qquad\text{in }Q_R(0).
\]
Moving along $\partial_x$, $Y$, $-\partial_x$ and $-Y$ over lengths $h$, $s$, $h$ and $s$ changes only the variable $y$, by $hs$. With $h=s=|\eta|^{1/2}$, the bounds on $\partial_xu_k$ and $Yu_k$ give $|u_k(x,y+\eta,t)-u_k(x,y,t)|\le C|\eta|^{1/2}$, and moving along $Y$ over the length $s$ and then correcting $y$ in the same way gives $|u_k(x,y,t+s)-u_k(x,y,t)|\le C(|s|+|xs|^{1/2})$, while $u_k$ is Lipschitz continuous in $x$. Hence the functions $u_k$ are uniformly Hölder continuous of order $\frac12$ on compact subsets of $Q_R(0)$. By the Arzelà--Ascoli theorem and a diagonal argument, a subsequence converges locally uniformly in $\R^3$ to a function $u_0$, and $\partial_{xx}u_k$ and $Yu_k$ converge weakly-$*$ in $L^\infty_{\mathrm{loc}}(\R^3)$ to $\partial_{xx}u_0$ and $Yu_0$.

The limit satisfies $0\le u_0\le\frac M2(x^+)^2$. Thus $u_0=0$ in $\{x\le0\}$ and $\sup_{Q_R(0)}u_0\le\frac M2R^2$ for $R\ge1$. Every compact subset of $\{x>0\}$ lies in $\{u_k>0\}$ for large $k$. There $a_k\to a(P)$ and $f_k\to f(P)$ uniformly, and $\rho_k\beta_k\,\partial_xu_k\to0$. Passing to the limit in \eqref{eq:rescaled}, we obtain $a(P)\,\partial_{xx}u_0+Yu_0=f(P)$ almost everywhere in $\{x>0\}$. By \cref{lem:liouville}, $u_0=\frac{f(P)}{2a(P)}(x^+)^2$. The limit does not depend on the subsequence, and therefore
\begin{equation}\label{eq:blowup}
 u_k\to\frac{f(P)}{2a(P)}(x^+)^2\qquad\text{locally uniformly in }\R^3.
\end{equation}

Taking $P_k=P$ shows that every point of $\Gamma\cap\overline{Q''}$ is regular, and by~(i) it is non-characteristic. Next we show that there are $c_1>0$ and $\rho_1\in(0,\rho_0]$ such that
\begin{equation}\label{eq:lower}
 u(x_P+\rho,y_P,t_P)\ge c_1\rho^2\qquad\text{for }P\in\Gamma\cap\overline{Q''}\text{ and }0<\rho<\rho_1.
\end{equation}
Otherwise there are $P_k\in\Gamma\cap\overline{Q''}$ and $\rho_k\searrow0$ with $u_k(1,0,0)\le k^{-1}$, since $\Lambda_k(1,0,0)=(x_{P_k}+\rho_k,y_{P_k},t_{P_k})$. This contradicts \eqref{eq:blowup} at a limit point $P$ of $(P_k)$, because $f(P)>0$.

By the compactness of $\Gamma\cap\overline{Q''}$ and the uniform size of the neighborhoods above, there is $\varepsilon>0$ with the following property. If $z=(x,y,t)\in\{u>0\}\cap Q'$ has distance less than $\varepsilon$ to $\Gamma$, then there is $P\in\Gamma\cap Q''$ with $y_P=y$, $t_P=t$ and $0<x-x_P<\rho_1$, namely $P=(\Phi(y,t),y,t)$ for one of the functions $\Phi$ above. Put $d=x-x_P$. Then $z=(x_P+d,y_P,t_P)\in Q_{\rho_0}(P)$, and \eqref{eq:quad} and \eqref{eq:lower} give $c_1d^2\le u(z)\le\frac M2d^2$. Since $\Gamma$ is a Lipschitz graph in the diffusive direction, $d$ is comparable to the distance from $z$ to $\Gamma$. This proves the comparability.
\end{proof}

\section{Euclidean regularity of the free boundary}\label{sec:euclid}
Throughout this section we keep the notation of the proof of \cref{thm:A}\,(ii). In particular, $\Gamma$ is the graph of a function $\Phi$ near each point of $\Gamma\cap\overline{Q''}$, and $M$, $\rho_0$ and $C_1$ are the constants in \eqref{eq:flatP} and \eqref{eq:quad}. The Euclidean regularity of $\Phi$ follows from boundary estimates for derivatives of $u$. We first show that these derivatives vanish linearly at $\Gamma$, and that $u_x$ vanishes at an exact linear rate.
\begin{lemma}\label{lem:linear}
There are $c_2>0$ and $\rho_2\in(0,\rho_0]$ with the following property. Let $(x,y,t)\in\{u>0\}$ be such that $(\Phi(y,t),y,t)\in\Gamma\cap\overline{Q''}$ and $x-\Phi(y,t)<\rho_2$. Then
\[
 c_2\bigl(x-\Phi(y,t)\bigr)\le u_x\le M\bigl(x-\Phi(y,t)\bigr)\qquad\text{and}\qquad|u_y|+|u_t|\le C_0M\bigl(x-\Phi(y,t)\bigr)
\]
at the point $(x,y,t)$.
\end{lemma}
\begin{proof}
The upper bounds follow from \eqref{eq:quad} and \eqref{eq:cone}. For the lower bound it suffices to find $c_2$ and $\rho_2$ such that $u_x(x_P+\rho,y_P,t_P)\ge c_2\rho$ for $P\in\Gamma\cap\overline{Q''}$ and $0<\rho<\rho_2$. Otherwise there are $P_k\in\Gamma\cap\overline{Q''}$ and $\rho_k\searrow0$ such that, in the notation of the proof of \cref{thm:A}\,(ii), $\partial_xu_k(1,0,0)\le k^{-1}$. We may assume that $P_k\to P$. For large $k$ the segment $[\frac12,\frac32]\times\{(0,0)\}$ lies in $\{u_k>0\}$, where $u_k$ is smooth and $|\partial_{xx}u_k|\le M$. Hence the function
\[
 g_k(s)=u_k(s,0,0)-\frac{f(P)}{2a(P)}s^2
\]
satisfies $|g_k''|\le M_1$ on $[\frac12,\frac32]$, where $M_1=M+f(P)/a(P)$. For $0<h\le\frac12$, Taylor's formula gives
\[
 |g_k'(1)|\le\frac2h\sup_{[\frac12,\frac32]}|g_k|+\frac{M_1h}{2}.
\]
By \eqref{eq:blowup}, the supremum tends to zero as $k\to\infty$. Since $h$ is arbitrary, $\partial_xu_k(1,0,0)\to f(P)/a(P)>0$, which is a contradiction.
\end{proof}

Next we show that the regularity of $\Phi$ follows from the behavior of two quotients of derivatives at $\Gamma$.
\begin{lemma}\label{lem:levelsets}
Let $P\in\Gamma\cap Q''$, and let $N=I\times D$ be a neighborhood of $P$, with an open interval $I$ and an open rectangle $D$, in which $\Gamma$ is the graph of $\Phi$ over $D$. Assume that there are functions $q_1$ and $q_2$ on $(\{u>0\}\cup\Gamma)\cap N$ with
\[
 |q_i(z)-q_i(z')|\le C_q|z-z'|^\gamma\qquad\text{for }i=1,2
\]
such that $u_y/u_x-q_1$ and $u_t/u_x-q_2$ tend to zero uniformly as $x-\Phi(y,t)\searrow0$. Then, after shrinking $N$, $\Phi\in C^{1,\gamma}(D)$ and $\nabla\Phi(y,t)=-\bigl(q_1,q_2\bigr)\bigl(\Phi(y,t),y,t\bigr)$.
\end{lemma}
\begin{proof}
After shrinking $N$, \cref{lem:linear} shows that $u_x>0$ in $\{u>0\}\cap N$ and that
\[
 \frac{c_2}{2}\bigl(x-\Phi(y,t)\bigr)^2\le u(x,y,t)\le\frac M2\bigl(x-\Phi(y,t)\bigr)^2\qquad\text{in }\{u>0\}\cap N.
\]
Hence, for small $\varepsilon>0$, the level set $\{u=\varepsilon\}\cap N$ is the graph of a function $\Phi_\varepsilon$ on $D$, and $0<\Phi_\varepsilon-\Phi\le(2\varepsilon/c_2)^{1/2}$. By the implicit function theorem, $\Phi_\varepsilon$ is smooth and
\[
 \nabla\Phi_\varepsilon(y,t)=-\Bigl(\frac{u_y}{u_x},\frac{u_t}{u_x}\Bigr)\bigl(\Phi_\varepsilon(y,t),y,t\bigr).
\]
As $\varepsilon\searrow0$, we have $\Phi_\varepsilon\to\Phi$ uniformly on $D$. Since $q_1$ and $q_2$ are uniformly continuous and the two quotients differ from them by a quantity that tends to zero uniformly as $x-\Phi(y,t)\searrow0$, $\nabla\Phi_\varepsilon$ converges uniformly to the function $-(q_1,q_2)(\Phi(y,t),y,t)$. Therefore $\Phi\in C^1(D)$ with the stated gradient. Since $\Phi$ is Lipschitz, the Hölder bound for $q_1$ and $q_2$ gives
\[
 |\nabla\Phi(y,t)-\nabla\Phi(y',t')|\le C\,|(y,t)-(y',t')|^\gamma.
\]
\end{proof}

The main step is a decay estimate for the oscillation of the quotients at $\Gamma$. The tool is the local comparison principle of Litsgård and Nyström~\cite[Lemma~9.1]{LN22}. We first state the results of~\cite{LN22} that are used in the proof of \cref{prop:C1gamma} below and in \cref{sec:asian}, in the special case of one diffusive variable, and we record the facts about the reference points that these proofs need. We use the coordinates $(x,y,s)$ with $s=-t$, in which the time is reversed, and the operator
\begin{equation}\label{eq:LNop}
 \mathcal{L}v=(a\,v_x)_x+x\,v_y-v_s,
\end{equation}
where $a\in C^\infty(\R^3)$ satisfies $\kappa^{-1}\le a\le\kappa$ for some $\kappa\ge1$ and $a=1$ outside a compact set. This is the operator of~\cite[(1.9)]{LN22} with $m=1$ and $A=a$, and the assumptions~\cite[(1.10)--(1.12)]{LN22} hold. In the proofs of \cref{prop:C1gamma,prop:asian-cone}, the operator $\mathcal{K}$ is brought into this form by an exponential conjugation and the time reversal, and its coefficient is modified away from the region under consideration so that these assumptions hold. Weak solutions of $\mathcal{L}v=0$ are understood as in~\cite{LN22}, and since $a$ is smooth, they are smooth by Hörmander's theorem.

The vector fields $\partial_x$ and $x\partial_y-\partial_s$ are invariant under the left translations $P'\mapsto P\circ P'$ with respect to the group law
\begin{equation}\label{eq:grouplaw}
 P\circ(x,y,s)=(x_P+x,\ y_P+y-s\,x_P,\ s_P+s)\qquad\text{for }P=(x_P,y_P,s_P).
\end{equation}
This is the group law of~\cite[(2.2)]{LN22}, and it is the time-reversed form of the translations in \eqref{eq:Lambda}. With the quasi-norm $\|(x,y,s)\|=|x|+|y|^{1/3}+|s|^{1/2}$, we use the symmetric quasi-distance
\[
 d(P,P')=\tfrac12\bigl(\|P'^{-1}\circ P\|+\|P^{-1}\circ P'\|\bigr)
\]
of~\cite[(2.8)]{LN22} and the balls $B_r(P)=\{P':\ d(P',P)<r\}$. In particular
\begin{equation}\label{eq:LNballs}
 d\bigl(P\circ(x,0,0),P\bigr)=|x|.
\end{equation}

A function $\psi\colon\R^2\to\R$ is a Lipschitz function with constant $M\ge1$ in the sense of~\cite[(2.12)]{LN22} if
\begin{equation}\label{eq:LNlip}
 |\psi(y,s)-\psi(y',s')|\le M\Bigl(|s-s'|^{\frac12}+\bigl|y-y'+(s-s')\psi(y',s')\bigr|^{\frac13}\Bigr)\qquad\text{for all }(y,s),(y',s')\in\R^2,
\end{equation}
and then $\Omega=\{(x,y,s):\ x>\psi(y,s)\}$ is a Lipschitz domain with constant $M$. For $P_0=(x_0,y_0,s_0)\in\partial\Omega$ and $r>0$, the sets of~\cite[(4.2)]{LN22} are
\begin{equation}\label{eq:LNsets}
\begin{aligned}
 Q_{M,r}&=\{(x,y,s):\ |x|<4Mr,\ |y|<r^3,\ |s|<r^2\},\qquad Q_{M,r}(P_0)=P_0\circ Q_{M,r},\\
 \Omega_r(P_0)&=Q_{M,r}(P_0)\cap\{\psi(y,s)<x<x_0+4Mr\},\qquad\Delta_r(P_0)=Q_{M,r}(P_0)\cap\{x=\psi(y,s)\},
\end{aligned}
\end{equation}
and by~\cite[(4.3)]{LN22} there is $c=c(M)\ge1$ such that $\Omega\cap B_{r/c}(P_0)\subset\Omega_r(P_0)\subset\Omega\cap B_{cr}(P_0)$. The sets $\Omega_r(P_0)$ decrease as $r$ decreases. In~\cite{LN22} the set $\Delta_r(P_0)$ is defined with $|x|<r$ in place of $|x|<4Mr$, so it is contained in the set defined here, which is the part of the graph in the boundary of $\Omega_r(P_0)$. \Cref{lem:LN,lem:single-scale} below assume boundary values on $\Delta_{2r}(P_0)$, and their hypotheses are therefore at least as strong as with the definition of~\cite{LN22}. The reference points of~\cite[(2.14)]{LN22} are
\begin{equation}\label{eq:refdef}
 \Xi^+_{\varrho}(P_0)=P_0\circ\Bigl(\Lambda\varrho,\ -\tfrac23\Lambda\varrho^3,\ \varrho^2\Bigr)\qquad\text{and}\qquad\Xi^-_{\varrho}(P_0)=P_0\circ\Bigl(\Lambda\varrho,\ \tfrac23\Lambda\varrho^3,\ -\varrho^2\Bigr),
\end{equation}
where $\Lambda=\Lambda(M)\ge1$ and $c_0=c_0(M)\ge1$ are the constants fixed in~\cite[Remark~4.2]{LN22}. By \eqref{eq:grouplaw} and \eqref{eq:LNlip}, the height of the graph below $\Xi^\pm_\varrho(P_0)$ satisfies
\[
 \bigl|\psi\bigl(y_0\mp\tfrac23\Lambda\varrho^3\mp\varrho^2x_0,\ s_0\pm\varrho^2\bigr)-x_0\bigr|\le M\Bigl(1+\bigl(\tfrac23\Lambda\bigr)^{\frac13}\Bigr)\varrho,
\]
and by~\cite[Lemma~4.12]{LN22} the reference points lie in $\Omega_r(P_0)$ at quasi-distance at least $\varrho/C$ from $\Delta_{2r}(P_0)$ for $0<\varrho<r/c_0$, where $C=C(M)\ge1$. By \eqref{eq:LNballs}, the segment $\{\Xi^\pm_\varrho(P_0)\circ(-x,0,0):\ 0\le x<\varrho/C\}$ lies in $Q_{M,r}(P_0)$ at quasi-distance less than $\varrho/C$ from $\Xi^\pm_\varrho(P_0)$, so it does not meet $\Delta_{2r}(P_0)$ and lies in $\Omega$. Hence the distance from $\Xi^\pm_\varrho(P_0)$ to the graph in the $x$-direction satisfies
\begin{equation}\label{eq:refpoints}
 \frac{\varrho}{C}\le x_{\Xi^\pm_\varrho(P_0)}-\psi\bigl(y_{\Xi^\pm_\varrho(P_0)},s_{\Xi^\pm_\varrho(P_0)}\bigr)\le C\varrho\qquad\text{for }0<\varrho<\frac{r}{c_0},
\end{equation}
after enlarging $C=C(M)$.

Following~\cite{Man97} and~\cite[Section~6]{LN22}, the propagation set of a point of an open set $D\subset\R^3$ consists of the points of $D$ that can be reached from it along curves in $D$ that follow $\pm\partial_x$ or the field $x\partial_y-\partial_s$, and the Kolmogorov boundary $\partial_KD$ is the set of points of $\partial D$ that lie in the closure of the propagation set of some point of $D$. For $D=\Omega_r(P_0)$ it consists of $\Delta_r(P_0)$, the face $x=x_0+4Mr$, the earliest time face, and the parts of the two $y$-faces of $Q_{M,r}(P_0)$ through which the field $x\partial_y-\partial_s$ leaves $D$, see~\cite[Section~6]{LN22}. The following lemma collects the solvability of the Dirichlet problem and the comparison principle that we use.
\begin{lemma}\label{lem:dirichlet}
Let $\Omega$ be a Lipschitz domain with constant $M$, let $P_0\in\partial\Omega$ and $r>0$, and write $D=\Omega_r(P_0)$.
\begin{enumerate}[label=\textup{(\roman*)}]
\item For every $\phi\in C(\overline D)$ there is a unique weak solution $Z\in C(\overline D)$ of $\mathcal{L}Z=0$ in $D$ with $Z=\phi$ on $\partial_KD$. It satisfies $\sup_D|Z|\le\sup_{\partial_KD}|\phi|$, and $Z\ge0$ if $\phi\ge0$ on $\partial_KD$.
\item Let $D'\subset\R^3$ be a bounded open set, and let $W\in C^2(D')\cap C(\overline{D'})$ satisfy $\mathcal{L}W\le0$ in $D'$ and $W\ge0$ on $\partial_KD'$. Then $W\ge0$ in $D'$.
\end{enumerate}
In particular, if $W\in C^2(D)\cap C(\overline D)$ satisfies $\mathcal{L}W\le0$ in $D$, then $W\ge Z$ in $D$ for the solution $Z$ of~\textup{(i)} with $\phi=W$, and the reverse inequality holds if $\mathcal{L}W\ge0$.
\end{lemma}
\begin{proof}
The existence, the uniqueness and the bound in~(i) are~\cite[Lemma~6.1]{LN22}, where $\phi$ is extended continuously to $\R^3$. If $\phi\ge0$, the bound applied to the data $\sup\phi-\phi$, whose solution is $\sup\phi-Z$, gives $Z\ge0$. Statement~(ii) is the weak maximum principle with respect to the Kolmogorov boundary, see~\cite{Man97}, \cite[Section~2.4]{NP16} and~\cite[Section~6]{LN22}. It holds for every bounded open set $D'$, because $s$ decreases strictly along the field $x\partial_y-\partial_s$, so that the propagation set of every point of $D'$ reaches $\partial D'$, and a supersolution that attains its minimum over $\overline{D'}$ at a point of $D'$ therefore attains it also on $\partial_KD'$ by the strong maximum principle of Bony~\cite{Bon69}. The last assertion follows from~(ii) applied to $W-Z$, which is smooth in $D$, continuous on $\overline D$ and zero on $\partial_KD$.
\end{proof}

The comparison principle of~\cite[Lemma~9.1]{LN22} reads as follows for $m=1$.
\begin{lemma}\label{lem:LN}
Let $\Omega$ be a Lipschitz domain with constant $M$. There is a constant $c=c(M,\kappa)\ge1$ with the following property. Let $P_0\in\partial\Omega$ and $r>0$, and let $u$ and $v$ be nonnegative weak solutions of $\mathcal{L}u=\mathcal{L}v=0$ in $\Omega_{2r}(P_0)$ that vanish continuously on $\Delta_{2r}(P_0)$, with $u>0$ in $\Omega_{2r}(P_0)$. Then
\[
 \frac1c\,\frac{v(\Xi^-_\varrho(P_0))}{u(\Xi^+_\varrho(P_0))}\le\frac{v(P)}{u(P)}\le c\,\frac{v(\Xi^+_\varrho(P_0))}{u(\Xi^-_\varrho(P_0))}\qquad\text{for }P\in\Omega_{\varrho/c}(P_0)\text{ and }0<\varrho\le\frac{r}{c_0^2}.
\]
\end{lemma}
\begin{proof}
This is~\cite[Lemma~9.1]{LN22} with $m=1$, in which $\varrho_1=r/c_0^2$. The positivity of $u$ is not listed among the hypotheses of~\cite[Lemma~9.1]{LN22}, where it is implicit in the quotients, and in \cref{prop:C1gamma,lem:single-scale} it is verified before the lemma is applied. The structural assumption~\cite[(3.4)]{LN22}, which requires $a$ and $\psi$ to be independent of $y$, is not needed here, as~\cite[Section~3.1]{LN22} states that ``up until and including Lemma 9.1, assumption (3.4) is not used''. Hence the lemma holds for the operator \eqref{eq:LNop} in every Lipschitz domain, with a constant that depends only on $M$ and $\kappa$.
\end{proof}

We now return to the quotients. In the boundary Harnack inequality of~\cite{Nys26a} the main difficulty is the balance between forward and backward reference values. For the quotients below this balance is automatic, because numerator and denominator are controlled by the distance to $\Gamma$ through \cref{lem:linear}. For the same reason, the terms that appear because $\partial_x$, $\partial_y$ and $\partial_t$ do not commute with $\mathcal{K}$, and the zero-order term that remains after the drift in $x$ is removed, can be treated as right-hand sides.
\begin{proposition}\label{prop:C1gamma}
Under the hypotheses of \cref{thm:A}, every point of $\Gamma\cap Q''$ has a neighborhood $N$ with the following property. There are functions $q_1$ and $q_2$ on $(\{u>0\}\cup\Gamma)\cap N$ that are Hölder continuous with respect to the Euclidean distance, such that $u_y/u_x-q_1$ and $u_t/u_x-q_2$ tend to zero uniformly as $x-\Phi(y,t)\searrow0$.
\end{proposition}
\begin{proof}
Fix $P\in\Gamma\cap Q''$ and write $d=x-\Phi(y,t)$. Differentiating $\mathcal{K}u=f$ in $\{u>0\}$ gives
\[
 \mathcal{K}u_x=f_x-a_xu_{xx}-\beta_xu_x-u_y\qquad\text{and}\qquad\mathcal{K}u_\ell=f_\ell-a_\ell u_{xx}-\beta_\ell u_x\quad\text{for }\ell\in\{y,t\}.
\]
In these identities we write $a_\ell u_{xx}=(a_\ell/a)\,au_{xx}$ for $\ell\in\{x,y,t\}$ and replace $au_{xx}$ by $f-\beta u_x-Yu+ru$. By \cref{lem:linear} and \eqref{eq:quad}, the functions $u$, $u_x$, $u_y$, $u_t$ and $Yu$ are bounded by multiples of $d$ near $P$. Hence, for $w=u_\ell$, the function $H_\ell=f_\ell-a_\ell f/a$ satisfies $|\mathcal{K}w-H_\ell|\le C_2d$ near $P$. Put $\widetilde w=w-H_\ell u/f$, and write $\widetilde{u_\ell}$ for this function when $w=u_\ell$. Since
\[
 \mathcal{K}\Bigl(\frac{H_\ell}{f}u\Bigr)=H_\ell+u\,(\mathcal{K}+r)\Bigl(\frac{H_\ell}{f}\Bigr)+2a\Bigl(\frac{H_\ell}{f}\Bigr)_xu_x,
\]
we obtain, after enlarging $C_2$, that $|\mathcal{K}\widetilde w|\le C_2d$ and $|\widetilde w-w|\le C_2d^2$ near $P$. In particular $\widetilde{u_x}$ is comparable to $d$. Since $u\le Md^2$ and $|u_\ell|\le Cd$, it follows that
\[
 \Bigl|\frac{u_\ell}{u_x}-\frac{\widetilde{u_\ell}}{\widetilde{u_x}}\Bigr|=\frac{u\,|H_\ell u_x-H_xu_\ell|}{f\,u_x\widetilde{u_x}}\le Cd\qquad\text{for }\ell\in\{y,t\}.
\]
It therefore suffices to find functions $q_1$ and $q_2$ on $N$ that are Hölder continuous with respect to the Euclidean distance, such that $\widetilde{u_y}/\widetilde{u_x}-q_1$ and $\widetilde{u_t}/\widetilde{u_x}-q_2$ tend to zero uniformly as $d\searrow0$. We obtain them by extending the boundary values of these two quotients on $\Gamma$.

Put $\widetilde\beta=\beta-a_x$, and let $\varphi$ be a smooth function near $P$ with $\varphi_x=-\widetilde\beta/(2a)$. A direct computation shows that
\[
 \mathcal{K}\bigl(e^{\varphi}v\bigr)=e^{\varphi}\bigl(\mathcal{L}v+cv\bigr)\qquad\text{with}\qquad\mathcal{L}v=(av_x)_x+Yv,
\]
where the function
\[
 c=a\bigl(\varphi_{xx}+\varphi_x^2\bigr)+\beta\varphi_x+Y\varphi-r
\]
is smooth. Put $V_0=e^{-\varphi}\widetilde{u_x}$, and let $v$ be $e^{-\varphi}\widetilde{u_y}$ or $e^{-\varphi}\widetilde{u_t}$. Then $V_0$ is comparable to $d$ and $|v|\le C_3d$. Since $cV_0$ and $cv$ are bounded by multiples of $d$ as well, we obtain $|\mathcal{L}V_0|+|\mathcal{L}v|\le C_3d$ after enlarging $C_3$. The function $q=v/V_0$ is the quotient in question. After the time reversal $s=-t$, the operator $\mathcal{L}$ takes the form \eqref{eq:LNop}. We may change $a$ away from a neighborhood of $P$ so that it stays smooth and uniformly elliptic and equals one outside a compact set, as required there. This does not affect the functions above near $P$.

Near the image $(x_0,y_0,s_0)$ of $P$, the set $\{u>0\}$ becomes $\{x>\widetilde\Phi(y,s)\}$ with $\widetilde\Phi(y,s)=\Phi(y,-s)$. The geometric part of the proof uses only the following property of $\widetilde\Phi$. There are $\delta\in(0,1]$ and $L\ge0$ such that
\begin{equation}\label{eq:halfLip}
 |\widetilde\Phi(y,s)-\widetilde\Phi(y',s')|\le L\bigl(|y-y'|+|s-s'|^{\frac12}\bigr)
\end{equation}
on the square $\Sigma$ of side $2\delta$ centered at $(y_0,s_0)$. Under the hypotheses of \cref{thm:A} this follows from~(i), because $|s-s'|\le2\delta\le2$ on $\Sigma$, so that $|s-s'|\le2|s-s'|^{1/2}$ there. Let $\pi$ be the nearest-point projection onto $\Sigma$, which clamps each coordinate to the corresponding interval, and put $\widehat\Phi=\widetilde\Phi\circ\pi$ on $\R^2$. Since $\pi$ does not increase the distance in each coordinate, $\widehat\Phi$ satisfies \eqref{eq:halfLip} on $\R^2$, and $|\widehat\Phi-x_0|\le m_0$ with $m_0=2L\delta^{1/2}$. It is also a Lipschitz function in the sense of \eqref{eq:LNlip}. Indeed, put
\[
 \Delta=\bigl|y-y'+(s-s')\widehat\Phi(y',s')\bigr|^{\frac13}+|s-s'|^{\frac12}\qquad\text{and}\qquad M_*=\max\bigl\{1,\ 2m_0,\ L(2+|x_0|+m_0)\bigr\}.
\]
If $\Delta\ge1$, then $|\widehat\Phi(y,s)-\widehat\Phi(y',s')|\le2m_0\le M_*\Delta$. If $\Delta<1$, both terms of $\Delta$ are smaller than one, and
\[
 |y-y'|\le\bigl|y-y'+(s-s')\widehat\Phi(y',s')\bigr|+(|x_0|+m_0)\,|s-s'|\le(1+|x_0|+m_0)\Delta.
\]
Together with \eqref{eq:halfLip}, this gives $|\widehat\Phi(y,s)-\widehat\Phi(y',s')|\le M_*\Delta$ again. Hence $\widehat\Omega=\{x>\widehat\Phi(y,s)\}$ is a Lipschitz domain with constant $M_*$ in the sense of \eqref{eq:LNlip}, and \cref{lem:dirichlet,lem:LN} apply to $\mathcal{L}$ in $\widehat\Omega$ with constants that depend only on $M_*$ and the ellipticity constant of $a$. We use the sets $\Omega_\rho(P')$ and the reference points $\Xi^\pm_\rho(P')$ of \eqref{eq:LNsets} and \eqref{eq:refdef} for graph points $P'$ near the image of $P$. For small $\rho$ these sets lie in the region where $\widehat\Phi=\widetilde\Phi$ and where the bounds above hold.

Fix such a point $P'$ and a small $\rho>0$, and let $c$ and $c_0$ be the constants of \cref{lem:LN}. Let $\underline q_\rho$ and $\overline q_\rho$ be the infimum and the supremum of $q$ on $\Omega_{2\rho}(P')$, and put $\omega_\rho=\overline q_\rho-\underline q_\rho$. The functions $v_1=v-\underline q_\rho V_0$ and $v_2=\overline q_\rho V_0-v$ are nonnegative in $\Omega_{2\rho}(P')$, vanish on the graph, and satisfy $|\mathcal{L}V_0|+|\mathcal{L}v_i|\le C_4d$, where $C_4$ does not depend on $\rho$ because $q$ is bounded. Let $\tau=s-s_*$, where $s_*$ is the smallest value of $s$ on $\Omega_{2\rho}(P')$, and put $\Psi=C_5\tau V_0$. Then $0\le\tau\le8\rho^2$, and
\[
 \mathcal{L}\Psi=C_5\bigl(\tau\,\mathcal{L}V_0-V_0\bigr)\le-C_4d\qquad\text{and}\qquad0\le\Psi\le C_6\rho^2V_0
\]
for a suitable $C_5$ and small $\rho$. Let $Z_1$, $Z_2$ and $Z_0$ be the solutions of $\mathcal{L}Z=0$ in $\Omega_{2\rho}(P')$ given by \cref{lem:dirichlet}\,(i) with the boundary values $v_1+\Psi$, $v_2+\Psi$ and $V_0$ on the Kolmogorov boundary. Since $v_i+\Psi$ and $V_0+\Psi$ are classical supersolutions, $V_0-\Psi$ is a classical subsolution, and $Z_1+Z_2-\omega_\rho Z_0$ has the boundary values $2\Psi\ge0$, \cref{lem:dirichlet} gives
\[
 Z_i\le v_i+\Psi,\qquad|Z_0-V_0|\le\Psi\qquad\text{and}\qquad Z_1+Z_2\ge\omega_\rho Z_0\qquad\text{in }\Omega_{2\rho}(P').
\]
In particular $Z_0$ is comparable to $d$, so that $Z_0>0$ in $\Omega_{2\rho}(P')$, as \cref{lem:LN} requires. Put $\rho_1=\rho/c_0^2$. By \eqref{eq:refpoints}, the reference points $\Xi^\pm_{\rho_1}(P')$ lie at a distance comparable to $\rho_1$ from the graph in the $x$-direction, so that $Z_0(\Xi^+_{\rho_1}(P'))\le C_7Z_0(\Xi^-_{\rho_1}(P'))$. One of the two functions, say $Z_1$, satisfies $Z_1(\Xi^-_{\rho_1}(P'))\ge\frac12\omega_\rho Z_0(\Xi^-_{\rho_1}(P'))$. Applying \cref{lem:LN} to $Z_1+\varepsilon Z_0$ and $Z_0$ and letting $\varepsilon\searrow0$, we obtain $Z_1\ge(2cC_7)^{-1}\omega_\rho Z_0$ in $\Omega_{\rho_1/c}(P')$. Hence
\[
 q-\underline q_\rho=\frac{v_1}{V_0}\ge\frac{Z_1-\Psi}{V_0}\ge\frac{\omega_\rho}{2cC_7}-C_8\rho^2\qquad\text{in }\Omega_{\rho_1/c}(P'),
\]
so that
\[
 \omega_{\lambda\rho}\le\Bigl(1-\frac{1}{2cC_7}\Bigr)\omega_\rho+C_8\rho^2\qquad\text{with}\qquad\lambda=\frac{1}{2cc_0^2}.
\]
Iterating along the geometric sequence of scales $\lambda^k\rho$ gives $\omega_\rho\le C\rho^{\gamma_0}$ for some $\gamma_0\in(0,1)$ and all small $\rho$, uniformly in $P'$.

It remains to construct $q_1$ and $q_2$ from the boundary values of $q$. Let $P'$ be a graph point as above. The sets $\Omega_{2\rho}(P')$ decrease as $\rho$ decreases, $q$ takes its values in $[\underline q_\rho,\overline q_\rho]$ on $\Omega_{2\rho}(P')$, and each of these sets contains the intersection of $\widehat\Omega$ with a neighborhood of $P'$. Since $\omega_\rho\to0$, the limit $q(P')$ of $q(z)$ as $z\in\widehat\Omega$ tends to $P'$ exists, and, after a change of the constant,
\begin{equation}\label{eq:trace}
 |q(z)-q(P')|\le C\rho^{\gamma_0}\qquad\text{for }z\in\Omega_\rho(P')
\end{equation}
and all small $\rho$, uniformly in $P'$. Let $P''$ be another graph point near the image of $P$, write $P'=(x',y',s')$ and $P''=(x'',y'',s'')$, and put
\[
 \varrho=|x''-x'|+\bigl|y''-y'+(s''-s')x'\bigr|^{\frac13}+|s''-s'|^{\frac12}.
\]
For $0<\xi<\varrho$, the point $z_\xi=(x''+\xi,y'',s'')$ lies in $\widehat\Omega$, and by the group law \eqref{eq:grouplaw} it is the translate by $P'$ of the point
\[
 \bigl(x''-x'+\xi,\ y''-y'+(s''-s')x',\ s''-s'\bigr),
\]
which lies in the box $Q_{M_*,2\varrho}$ of \eqref{eq:LNsets}. Hence $z_\xi\in\Omega_{2\varrho}(P')$, and \eqref{eq:trace} gives $|q(z_\xi)-q(P')|\le C\varrho^{\gamma_0}$ when $\varrho$ is small. Letting $\xi\searrow0$, we obtain $|q(P'')-q(P')|\le C\varrho^{\gamma_0}$ when $\varrho$ is small, and since $q$ is bounded, this bound holds for all such pairs of graph points. Since $\varrho\le C|P''-P'|^{1/3}$ near the image of $P$, the boundary values of $q$ are Hölder continuous with exponent $\gamma_0/3$ with respect to the Euclidean distance.

Let $p_1$ and $p_2$ be the boundary values obtained in this way from $\widetilde{u_y}/\widetilde{u_x}$ and $\widetilde{u_t}/\widetilde{u_x}$, written in the original variables, and put
\[
 q_i(x,y,t)=p_i\bigl(\Phi(y,t),y,t\bigr)\qquad\text{for }i=1,2.
\]
Since $\Phi$ is Lipschitz, $q_1$ and $q_2$ are Hölder continuous with exponent $\gamma_0/3$ on a neighborhood $N$ of $P$ of the form in \cref{lem:levelsets}. Let $z=(x,y,t)\in\{u>0\}\cap N$ and $d=x-\Phi(y,t)$. In the time-reversed variables, $z$ is the translate by the image $P'$ of $(\Phi(y,t),y,t)$ of the point $(d,0,0)$, which lies in $Q_{M_*,2d}$. Hence $z\in\Omega_{2d}(P')$, and \eqref{eq:trace} together with the bound for $u_\ell/u_x-\widetilde{u_\ell}/\widetilde{u_x}$ gives
\[
 \Bigl|\frac{u_y}{u_x}-q_1\Bigr|+\Bigl|\frac{u_t}{u_x}-q_2\Bigr|\le C\bigl(d^{\gamma_0}+d\bigr)\qquad\text{at }z.
\]
The right-hand side tends to zero uniformly as $d\searrow0$, which proves the proposition.
\end{proof}

\begin{proof}[Proof of \cref{thm:A}\,(iii)]
By \cref{prop:C1gamma} and its proof, the functions $q_1$ and $q_2$ are Hölder continuous with the exponent $\gamma_0/3$, where $\gamma_0$ is the decay exponent of the oscillation of the quotients. Hence \cref{lem:levelsets} shows that $\Phi$ is of class $C^{1,\gamma}$ with $\gamma=\gamma_0/3$ near every point of $\Gamma\cap Q''$, and $\Phi_y$ and $\Phi_t$ are the limits of $-u_y/u_x$ and $-u_t/u_x$ at $\Gamma$. Since $\Gamma\cap\overline{Q'}$ is a compact subset of $Q''$, the exponent and the constants can be chosen uniformly.
\end{proof}

The method of this section does not go beyond Hölder continuous first derivatives.
\begin{remark}\label{rem:higher}
We expect that $\Phi$ is $C^\infty$ when the data are smooth. By \cref{lem:levelsets}, this would follow from a higher-order boundary Harnack inequality, that is, from the $C^{k,\gamma}$ regularity of $u_y/u_x$ and $u_t/u_x$ up to $\Gamma$ whenever $\Phi\in C^{k,\gamma}$, as in the parabolic result of Kukuljan~\cite[Theorem~1.3]{Kuk22}. The partial hodograph transform $(x,y,t)\mapsto(u_x,y,t)$ flattens $\Gamma$, but then the transport velocity depends on the solution. The change of variables $x\mapsto x-\Phi(y,t)$ also flattens $\Gamma$, and the transport velocity becomes $x+\Phi(y,t)$. In both forms one needs a Schauder estimate up to a boundary that is a graph in the diffusive variable, and we are not aware of such an estimate.
\end{remark}

\section{The American Asian option}\label{sec:asian}
Throughout this section $x=S$ and $y=I$, so that the transport term of \eqref{eq:K} is $S\partial_I$. We write points as $(t,S,I)$, in the order that is customary in finance, so that the diffusive variable comes second. In the region $\{A>K\}$ the obstacle is $\psi=I/t-K$, which does not depend on $S$. We first identify the normalized problem and the region where exercise can occur.
\begin{lemma}\label{lem:asian-incentive}
In $\{A>K\}$ the function $u=V-\psi$ satisfies $u\ge0$ and $\mathcal{K}u=f\chi_{\{u>0\}}$ almost everywhere, with
\[
 a=\frac{\sigma^2S^2}{2},\qquad\beta=(r-q)S,\qquad f=\frac{A-S}{t}+r(A-K).
\]
In particular, exercise occurs only where $S\le A+rt(A-K)$.
\end{lemma}
\begin{proof}
By~\cite[Theorem~1.1]{MP09}, the value function is a strong solution of the obstacle problem, with $V_S$, $V_{SS}$ and $YV$ in $L^p_{\mathrm{loc}}$ for every $p<\infty$, and by the optimal regularity theorem~\cite[Theorem~1.7]{FNPP10}, applied in $\{A>K\}$, where the obstacle is smooth, $V_{SS}$ and $YV$ are locally bounded in $\{0<t<T,\ S>0,\ A>K\}$. These results are stated for $q=0$. The dividend yield only changes the coefficient of $\partial_S$, and the class of operators of~\cite{DFPP08,FNPP10}, on which~\cite{MP09} rests, contains first-order terms with Hölder continuous coefficients, so the results apply to the operator with the coefficients above. That is, $\mathcal{K}V\le0$ and $V\ge\psi$ almost everywhere in $\{A>K\}$, and $\mathcal{K}V=0$ in $\{V>\psi\}$. In particular $V_S$ is continuous. Since $\psi_t=-I/t^2$, $\psi_I=1/t$, and $\psi$ does not depend on $S$, we obtain
\[
 \mathcal{K}\psi=-\frac{I}{t^2}+\frac{S}{t}-r\Bigl(\frac{I}{t}-K\Bigr)=\frac{S-A}{t}-r(A-K)=-f.
\]
Hence $\mathcal{K}u=f$ in $\{u>0\}$. In the contact set $\mathcal{K}u=0$ almost everywhere, and $\mathcal{K}\psi=\mathcal{K}V\le0$ there, so $f\ge0$ almost everywhere on the contact set. This holds at every point of the contact set. Indeed, if $u(P)=0$ and $f(P)<0$, the contact set has measure zero near $P$, so that $\mathcal{K}u=f$ near $P$ in the sense of distributions, and $u$ is smooth near $P$ by Hörmander's theorem~\cite{Hor67}. Since $u$ attains its minimum $0$ at $P$, this gives $\mathcal{K}u(P)=au_{SS}(P)\ge0>f(P)$, a contradiction.
\end{proof}

The scaling of the problem gives an identity that reduces the cone estimate in the $I$-direction to an estimate in the strike.
\begin{lemma}\label{lem:euler}
For every $\lambda>0$ we have $V(t,\lambda S,\lambda I;\lambda K)=\lambda V(t,S,I;K)$. Consequently
\[
 SV_S+IV_I+KV_K=V
\]
wherever $V$ is differentiable, and $u=V-\psi$ satisfies
\[
 -I\,u_I=S\,u_S-u+K\,\partial_Ku.
\]
\end{lemma}
\begin{proof}
The processes started from $(\lambda S,\lambda I)$ are $(\lambda S_s,\lambda I_s)$, and $(\lambda I_\theta/\theta-\lambda K)^+=\lambda(I_\theta/\theta-K)^+$ for every stopping time $\theta$. This gives the homogeneity, and differentiation at $\lambda=1$ gives Euler's identity. Since $u_S=V_S$, $u_I=V_I-1/t$ and $\partial_Ku=V_K+1$, we obtain
\[
 -I\,u_I=\frac{I}{t}-IV_I=\frac{I}{t}-V+SV_S+KV_K=S\,u_S-u+K\,\partial_Ku.
\]
\end{proof}

The value function inherits monotonicity from the payoff. Together with \cref{lem:euler}, the monotonicity in the strike becomes a dilation property of the exercise region.
\begin{lemma}\label{lem:asian-monotone}
For $h\ge0$ the value function satisfies
\[
 V(t,S+h,I)\ge V(t,S,I)\qquad\text{and}\qquad V(t,S,I)\le V(t,S,I+h)\le V(t,S,I)+\frac ht,
\]
and $V(t,S,I;K)-h\le V(t,S,I;K+h)\le V(t,S,I;K)$. Moreover, $V$ and $u$ are convex in $S$. Let $\mathcal{E}$ be the exercise region in $\{A>K\}$. If $(t,S,I)\in\mathcal{E}$, then $(t,S',I')\in\mathcal{E}$ whenever $S'\le S$ and $I'\ge I$, and $(t,\lambda S,\lambda I)\in\mathcal{E}$ for every $\lambda\ge1$.
\end{lemma}
\begin{proof}
Write $I_\theta=I+SJ_\theta$, where $J_\theta=\int_t^\theta S_\tau/S\,d\tau$ depends neither on $S$ nor on $I$. For every stopping time $\theta$ the payoff $(I_\theta/\theta-K)^+$ is nondecreasing in $S$ and $I$ and nonincreasing in $K$. Since $I_\theta$ is affine in $S$, the payoff is also convex in $S$. It increases by at most $h/\theta\le h/t$ when $I$ is replaced by $I+h$, and it decreases by at most $h$ when $K$ is replaced by $K+h$. Since the discount factor is at most one, taking expectations and suprema gives the inequalities for $V$ and preserves the convexity in $S$. Since $\psi$ does not depend on $S$, the function $u$ is convex in $S$ as well.

In $\{A>K\}$ we have $u=V-I/t+K$. Let $(t,S,I)\in\mathcal{E}$, so that $u(t,S,I)=0$. Then $u(t,S',I)\le u(t,S,I)=0$ for $S'\le S$, and
\[
 u(t,S,I+h)\le V(t,S,I)+\frac ht-\frac{I+h}{t}+K=u(t,S,I)=0.
\]
Since $u\ge0$, the points $(t,S',I)$ and $(t,S,I+h)$ lie in $\mathcal{E}$, and the first assertion follows. Write $u^K$ and $\mathcal{E}^K$ to indicate the strike. The bounds in $K$ give $u^{K+h}\ge u^K$, so $\mathcal{E}^{K+h}\subset\mathcal{E}^K$. By \cref{lem:euler}, $(t,S,I)\in\mathcal{E}^K$ if and only if $(t,\lambda S,\lambda I)\in\mathcal{E}^{\lambda K}$. For $\lambda\ge1$ we have $\mathcal{E}^{\lambda K}\subset\mathcal{E}^K$, which proves the last assertion.
\end{proof}

Consequently $u$ is diffusive-monotone, and the exercise region in $\{A>K\}$ is $\{S\le\Phi(I,t)\}$ for a function $\Phi$ with $\Phi(I+h,t)\ge\Phi(I,t)$ and $\Phi(\lambda I,t)\ge\lambda\Phi(I,t)$ for $h\ge0$ and $\lambda\ge1$. In particular $u_I\le0\le u_S$ wherever these derivatives exist. By \cref{lem:asian-gap,lem:asian-boundary} below, $\Phi$ is positive and continuous over compact subsets of $\{0<t<T,\ A>K\}$, and the free boundary there is the graph $S=\Phi(I,t)$.

Two further comparisons of payoffs bound $-u_I$ and $u_t$ by multiples of $Su_S-u$. The first one raises the average and the strike by the same amount, and the second one moves the starting time forward.
\begin{lemma}\label{lem:asian-shift}
Let $h\ge0$ with $t+h\le T$. Then
\[
 V(t,S,I+th;K+h)\le V(t,S,I;K)\qquad\text{and}\qquad V\Bigl(t+h,S,\frac{(t+h)^2}{t^2}\,I;K+\frac{hI}{t^2}\Bigr)\le V(t,S,I;K).
\]
Consequently, in $\{u>0\}\cap\{A>K\}$,
\[
 0\le-t(A-K)\,u_I\le Su_S-u\qquad\text{and}\qquad u_t\le\frac{A}{t\min\{K,A-K\}}\bigl(Su_S-u\bigr).
\]
\end{lemma}
\begin{proof}
Let $J_\theta$ be as in the proof of \cref{lem:asian-monotone}. For a stopping time $\theta\ge t$, replacing $I$ by $I+th$ and $K$ by $K+h$ changes the payoff into
\[
 \Bigl(\frac{I+SJ_\theta}{\theta}-K-h\Bigl(1-\frac t\theta\Bigr)\Bigr)^+\le\Bigl(\frac{I+SJ_\theta}{\theta}-K\Bigr)^+,
\]
which gives the first inequality. Put $I'=(t+h)^2I/t^2$ and $K'=K+hI/t^2$. Since the stock price is time homogeneous, $V(t+h,S,I';K')$ is the supremum of
\[
 \mathbb{E}\Bigl[e^{-r(\theta-t)}\Bigl(\frac{I'+SJ_\theta}{\theta+h}-K'\Bigr)^+\Bigr]
\]
over stopping times $\theta$ with values in $[t,T-h]$. The function
\[
 g(\theta)=\frac{(t+h)^2}{\theta+h}-\frac{t^2}{\theta}
\]
satisfies $g(t)=h$ and $g'\le0$ on $[t,\infty)$, because $t/\theta\le(t+h)/(\theta+h)$ there. Since $J_\theta\ge0$, we obtain
\[
 \frac{I'+SJ_\theta}{\theta+h}-K'\le\frac{I}{t^2}\bigl(g(\theta)-h\bigr)+\frac{I+SJ_\theta}{\theta}-K\le\frac{I+SJ_\theta}{\theta}-K.
\]
The supremum that defines $V(t,S,I;K)$ is taken over the larger set of stopping times with values in $[t,T]$, which gives the second inequality.

In $\{A>K\}$ we have $u=V-I/t+K$, and both shifts leave $I/t-K$ unchanged. Hence the two inequalities also hold for $u$. In $\{u>0\}$ the function $u$ is smooth, also in $K$, because $u(t,S,I;K)=K\,u(t,S/K,I/K;1)$ by \cref{lem:euler}, and differentiating at $h=0$ gives
\[
 t\,u_I+\partial_Ku\le0\qquad\text{and}\qquad u_t+2A\,u_I+\frac At\,\partial_Ku\le0.
\]
By \cref{lem:euler}, $K\partial_Ku=u-Su_S-Iu_I$. Inserting this into the first inequality and using $u_I\le0$ gives the bound for $u_I$. Inserting it into the second one gives
\[
 u_t\le\frac{A}{tK}\bigl(Su_S-u\bigr)+A\Bigl(\frac AK-2\Bigr)u_I.
\]
If $A\ge2K$, the last term is nonpositive. If $A<2K$, the bound for $u_I$ shows that it is at most
\[
 \frac{A(2K-A)}{tK(A-K)}\,\bigl(Su_S-u\bigr).
\]
In both cases the bound for $u_t$ follows.
\end{proof}

The next lemma shows that the exercise region stays inside the region where $f$ is positive, and that it contains a neighborhood of $S=0$.
\begin{lemma}\label{lem:asian-gap}
Let $E\subset\{(t,A):\ 0<t<T,\ A>K\}$ be compact. Then the following statements hold.
\begin{enumerate}[label=\textup{(\roman*)}]
\item The exercise region lies in $\{S<A+rt(A-K)\}$. In particular $f$ has a positive lower bound on the part of the exercise region over $E$.
\item There is $s_0>0$ such that $(t,S,tA)$ lies in the exercise region whenever $(t,A)\in E$ and $0<S\le s_0$.
\end{enumerate}
\end{lemma}
\begin{proof}
By \cref{lem:asian-incentive}, the exercise region lies in $\{f\ge0\}$. Hence the open set $D=\{f<0\}\cap\{A>K\}$ is contained in $\{u>0\}$, where $u$ is smooth and $\mathcal{K}u=f<0$. Suppose that $P$ lies in the exercise region with $f(P)=0$ and $0<t_P<T$. Since $f_S=-1/t\neq0$, the boundary of $D$ is smooth near $P$, so there is a ball $B=B_R(C)$ with $\overline B\subset D\cup\{P\}$ and $P\in\partial B$. Since $D$ lies on the side of larger $S$, the $S$-component $\nu_S$ of $\nu=(C-P)/R$ is positive. Put
\[
 w(z)=e^{-\alpha|z-C|^2}-e^{-\alpha R^2}.
\]
For small $\varepsilon>0$ we have $|S-C_S|\ge\frac12R\nu_S$ in $U=B\cap B_\varepsilon(P)$. Then the term $4\alpha^2a(S-C_S)^2e^{-\alpha|z-C|^2}$ in $\mathcal{K}w$ dominates the remaining terms for large $\alpha$, so that $\mathcal{K}w>0$ in $U$. Let $m>0$ be the minimum of $u$ on the compact set $\partial B_\varepsilon(P)\cap\overline B$. Then $u-mw\ge0$ on $\partial U$ and $\mathcal{K}(u-mw)<0$ in $U$. At a negative interior minimum we would have $\mathcal{K}(u-mw)\ge-r(u-mw)\ge0$, so $u\ge mw$ in $U$. Since $\partial_Sw(P)=2\alpha R\nu_Se^{-\alpha R^2}>0$, this gives $u(t_P,S_P+h,I_P)\ge c\,h$ for some $c>0$ and all small $h>0$. On the other hand, by \cref{lem:asian-monotone} the function $u$ vanishes on $\{(t_P,S,I_P):\ S\le S_P\}$, so that $u_S(P)=0$ by the continuity of $u_S$. Since $u_{SS}$ is bounded near $P$ by some $M$, we obtain $u(t_P,S_P+h,I_P)\le\frac M2h^2$. This is a contradiction for small $h$. The closure of the part of the exercise region over $E$ is compact, and by what we have shown $f$ is positive on it, also at $S=0$, where $f=A/t+r(A-K)$. This proves (i).

For (ii), let $J_\theta$ be as in the proof of \cref{lem:asian-monotone}. Since $(I_\theta/\theta-K)^+\le(I/t-K)^++SJ_\theta/t$ for every stopping time $\theta$ and $J_\theta\le J_T$, we have $0\le u\le C_{11}S$ in $\{A>K\}\cap\{t\ge t_*\}$, where $t_*>0$ is the smallest value of $t$ on $E$ and $C_{11}=Te^{|r-q|T}/t_*$. Fix $(t_0,A_0)\in E$ and $\delta>0$ with $t_0+\delta<T$ and $A_0-\delta>K$. For $0<s_1<(A_0-\delta)/2$ consider
\[
 N=\bigl\{(t,S,I):\ t_0<t<t_0+\delta,\ 0<S<s_1,\ A_0-\delta<A<A_0+\delta\bigr\},
\]
where $f\ge f_1>0$ with $f_1$ independent of $s_1$. Put
\[
 \overline u=\frac{k_1}{s_1^3}\Bigl(\Bigl(S-\frac{s_1}2\Bigr)^+\Bigr)^3+k_2(t-t_0)+\frac{k_3}{\delta}\bigl((A_0-A)^+\bigr)^2
\]
with $k_1=8C_{11}s_1$, $k_2=f_1/3$ and $k_3=f_1t_0/(6(A_0+\delta))$. This function is twice continuously differentiable in $S$ and continuously differentiable in $t$ and $I$. Since $YA=(S-A)/t$, $0<(A-S)/t<(A_0+\delta)/t_0$ and $(A_0-A)^+<\delta$ in $N$, we have
\[
 \mathcal{K}\overline u\le k_1\Bigl(\frac32\sigma^2+|r-q|\Bigr)+k_2+\frac{2k_3(A_0+\delta)}{t_0}\le f\qquad\text{in }N
\]
when $s_1$ is small. Where $S=s_1$, $t=t_0+\delta$ or $A=A_0-\delta$, we have $\overline u\ge C_{11}s_1\ge u$ when $s_1$ is small, and $u=0\le\overline u$ where $S=0$. We claim that $u\le\overline u$ in $N$. Otherwise, for small $\varepsilon>0$, the function $w_\varepsilon=u-\overline u-\varepsilon(t_0+\delta-t)$ has a positive maximum over $\overline N$ at some point $P$. Then $u(P)>0$, so $u$ is smooth near $P$, while $\overline u$ is twice continuously differentiable in $S$ and continuously differentiable in $t$ and $I$. The point $P$ does not lie on the faces $S=0$, $S=s_1$, $t=t_0+\delta$ and $A=A_0-\delta$, where $w_\varepsilon\le0$. Since $\mathcal{K}$ contains $+\partial_t$ and $YA=(S-A)/t<0$ in $N$, the transport field $Y=\partial_t+S\partial_I$ points into $N$ at every other point of $\overline N$, and the $S$-direction is tangential to the faces $t=t_0$ and $A=A_0+\delta$. Hence $Yw_\varepsilon(P)\le0$, $\partial_Sw_\varepsilon(P)=0$ and $\partial_{SS}w_\varepsilon(P)\le0$, so that $\mathcal{K}w_\varepsilon(P)\le-rw_\varepsilon(P)\le0$. On the other hand $\mathcal{K}w_\varepsilon=f-\mathcal{K}\overline u+\varepsilon+r\varepsilon(t_0+\delta-t)>0$ at $P$, a contradiction. In the terminology of \cref{sec:euclid}, the faces $S=0$, $S=s_1$, $t=t_0+\delta$ and $A=A_0-\delta$ form the Kolmogorov boundary of $N$, through which $Y$ leaves $N$. Hence $u\le\overline u$ in $N$, and in particular $u(t_0,S,t_0A_0)=0$ for $S\le s_1/2$. The constants can be chosen uniformly for $(t_0,A_0)\in E$, which proves (ii) with $s_0=s_1/2$.
\end{proof}

The bounds of \cref{lem:asian-shift} control the free boundary in $I$ and its descent in $t$. A nondegeneracy argument then controls its ascent in $t$. The ascent bound is only of order $h^{1/2}$, because \cref{lem:asian-shift} gives no lower bound for $u_t$, and a free boundary that moves up in $t$ can only be controlled through the nondegeneracy of $u$. This weaker bound is all that \cref{prop:asian-cone} needs, since it yields \eqref{eq:halfLip}.
\begin{lemma}\label{lem:asian-boundary}
Let $E\subset\{(t,A):\ 0<t<T,\ A>K\}$ be compact. There are constants $L\ge1$ and $h_0,\eta_0,c_3>0$ such that
\[
 0\le\Phi(I+h,t)-\Phi(I,t)\le Lh,\qquad-Lh\le\Phi(I,t+h)-\Phi(I,t)\le Lh^{\frac12}
\]
and
\[
 u\bigl(t,\Phi(I,t)+\eta,I\bigr)\ge c_3\eta^2
\]
whenever $(t,I/t)\in E$, $0\le h\le h_0$ and $0<\eta\le\eta_0$.
\end{lemma}
\begin{proof}
Let $E'$ be a compact neighborhood of $E$ in $\{0<t<T,\ A>K\}$. Let $t_*$ and $A_*$ be the smallest values of $t$ and $A$ on $E'$, and let $A^*$ be the largest value of $A$. By \cref{lem:asian-incentive,lem:asian-gap}, there are $s_0>0$ and $S_1$ such that $s_0\le\Phi(I,t)\le S_1$ whenever $(t,I/t)\in E'$. Let $N$ be the compact set of points $(t,S,I)$ with $(t,I/t)\in E'$ and $s_0/2\le S\le S_1+1$, and let $M$ bound $|u_{SS}|$ and $|Yu|$ on $N$. By \cref{lem:asian-shift},
\[
 -u_I\le L_1u_S\qquad\text{and}\qquad u_t\le L_1u_S\qquad\text{in }N\cap\{u>0\}
\]
with
\[
 L_1=\frac{(S_1+1)(1+A^*)}{t_*\min\{K,A_*-K\}}.
\]
By \cref{lem:asian-monotone}, $|u_I|\le1/t$, so that $|u_t|\le|Yu|+S|u_I|\le M_2$ in $N$ with $M_2=M+(S_1+1)/t_*$. Since $u_S$, $u_I$ and $u_t$ are bounded on $N$, the function $u$ is locally Lipschitz there. By \cref{lem:asian-gap}\,(i) and the uniform continuity of $f$ on $N$, there are $f_0>0$ and $\eta_1>0$ such that $f\ge f_0$ at every point of $N$ whose distance to the exercise region over $E'$ is less than $\eta_1$. Let $h_1>0$ be so small that all points used below for $0\le h\le h_1$ lie in $N$.

Let $(t,I/t)\in E$ and $0\le h\le h_1$. In $N\cap\{u>0\}$, the derivative of $u$ in the direction of $(0,L_1,1)$ in the variables $(t,S,I)$ is $L_1u_S+u_I\ge0$. As in the proof of \cref{thm:A}\,(i), $u$ is therefore nondecreasing along segments in this direction that start in $\{u>0\}$. Applying this to the segment from $(t,S,I)$ to $(t,S+L_1h,I+h)$ with $S>\Phi(I,t)$ and letting $S\searrow\Phi(I,t)$, we obtain $\Phi(I+h,t)\le\Phi(I,t)+L_1h$. The lower bound $\Phi(I+h,t)\ge\Phi(I,t)$ is part of \cref{lem:asian-monotone}. Similarly, the derivative in the direction of $(1,-L_1,0)$ is $u_t-L_1u_S\le0$ in $N\cap\{u>0\}$. If $u$ were positive at some point of the segment from $(t,\Phi(I,t),I)$ to $(t+h,\Phi(I,t)-L_1h,I)$, then, starting from the last zero of $u$ before that point, $u$ would be positive and nonincreasing along the segment, which is impossible. Hence $u$ vanishes on this segment, and $\Phi(I,t+h)\ge\Phi(I,t)-L_1h$.

Next we prove the lower bound for $u$. Let $(t_0,I_0/t_0)\in E$, put $\Phi_0=\Phi(I_0,t_0)$, and let $0<\eta\le\eta_0$, where $\eta_0\le4h_1^{1/2}$ is chosen small enough for the conditions below. Put $\rho=\eta/4$ and $S_0=\Phi_0+\eta/2$, and consider the half-cylinder
\[
 \mathcal{C}=\bigl\{(t,S,I):\ t_0<t<t_0+\rho^2,\ |S-S_0|<\rho,\ |I-I_0-S_0(t-t_0)|<\rho^3\bigr\}
\]
and the function
\[
 p=k_4(S-S_0)^2+k_5(t-t_0)+\frac{k_6}{\rho^4}\bigl(I-I_0-S_0(t-t_0)\bigr)^2,
\]
where $k_4=f_0/(8a_1)$, $k_5=f_0/4$, $k_6=f_0/16$ and $a_1$ is the maximum of $a$ on $N$. Since $Y(I-I_0-S_0(t-t_0))=S-S_0$, we have $\mathcal{K}p\le\frac58f_0+2k_4|\beta|\rho<f_0$ in $\mathcal{C}$ for small $\eta_0$. Moreover $p\ge c_4\rho^2$ on $\partial \mathcal{C}\setminus F$, where $F=\{t=t_0,\ |S-S_0|<\rho,\ |I-I_0|<\rho^3\}$ is the bottom face and $c_4=\min\{k_4,k_5,k_6\}$. For small $\eta_0$ the closure of $\mathcal{C}$ lies within distance $\eta_1$ of $(t_0,\Phi_0,I_0)$, so that $\mathcal{K}(u-p)>0$ in $\mathcal{C}\cap\{u>0\}$. The function $u-p$ is positive at $(t_0,S_0,I_0)$ and nonpositive on $\Gamma$. Its maximum over the closure of $\mathcal{C}\cap\{u>0\}$ is therefore positive, and it is attained at a point $(t_1,S_2,I_2)$ of $\partial \mathcal{C}\setminus F$. Indeed, at a positive maximum point that lies in $\mathcal{C}\cup F$ or on a part of $\partial \mathcal{C}$ through which $Y$ enters $\mathcal{C}$, the function $u$ is smooth, $Y(u-p)\le0$, $\partial_S(u-p)=0$ and $\partial_{SS}(u-p)\le0$, so that $\mathcal{K}(u-p)\le0$, which is impossible. Then
\[
 u(t_1,S_2,I_2)\ge c_4\rho^2,\qquad S_2<\Phi_0+\eta,\qquad t_1\le t_0+\rho^2\qquad\text{and}\qquad I_2\ge I_0-\rho^3.
\]
Since $u$ is nondecreasing in $S$, nonincreasing in $I$ and $|u_I|\le1/t$, this gives $u(t_1,\Phi_0+\eta,I_0)\ge c_4\rho^2-\rho^3/t_*$. For $t_0\le t\le t_1$, the bound for the descent of $\Phi$ shows that $\Phi_0+\eta-\Phi(I_0,t)\le\eta+L_1\rho^2$. The function $t\mapsto u(t,\Phi_0+\eta,I_0)$ is Lipschitz. Since $u_S$ vanishes on $\Gamma$ and $|u_{SS}|\le M$, its derivative is at most $L_1u_S\le L_1M(\eta+L_1\rho^2)$ where $u>0$, and it vanishes almost everywhere on the set where $u=0$. Hence
\[
 u(t_0,\Phi_0+\eta,I_0)\ge c_4\rho^2-\frac{\rho^3}{t_*}-L_1M\bigl(\eta+L_1\rho^2\bigr)\rho^2\ge\frac{c_4}2\rho^2=\frac{c_4}{32}\eta^2
\]
for small $\eta_0$, and the lower bound for $u$ holds with $c_3=c_4/32$.

Finally, let $h_0=\min\{h_1,c_3\eta_0^2/(2M_2)\}$ and $0<h\le h_0$, and put $\Delta=\Phi(I,t+h)-\Phi(I,t)$. If $\Delta>0$, let $\eta=\min\{\Delta/2,\eta_0\}$. Then $u(t+h,\Phi(I,t)+\eta,I)=0$, and the lower bound for $u$ together with $|u_t|\le M_2$ gives $c_3\eta^2\le M_2h\le c_3\eta_0^2/2$. Hence $\eta=\Delta/2$ and $\Delta\le2(M_2h/c_3)^{1/2}$. The lemma follows with $L=\max\{1,L_1,2(M_2/c_3)^{1/2}\}$.
\end{proof}

The cone estimate is the main point of the application. It is a local statement in the variables $t$ and $A$.
The lower bound for $u_t$ comes from a comparison at a single scale, which we isolate in the setting of \eqref{eq:LNop}. The constants in its conclusion depend neither on the function that is compared nor on the neighborhood of the free boundary in which it is defined, which is what the limit $h\searrow0$ in the proof of \cref{prop:asian-cone} requires.
\begin{lemma}\label{lem:single-scale}
Let $\Omega$ be a Lipschitz domain with constant $M$, let $\mathcal{L}$ be as in \eqref{eq:LNop}, and let $c$ and $c_0$ be the constants of \cref{lem:LN}. Let $P'\in\partial\Omega$ and $\rho>0$, and put $D=\Omega_{2\rho}(P')$, $\rho_1=\rho/c_0^2$ and $d=x-\psi(y,s)$. Let $U$ be an open neighborhood of $\overline D$. Let $w_0\in C(\overline D)$ be smooth in $D$, and let $w$ be Lipschitz continuous in $\Omega\cap U$ with $w_{xx}\in L^\infty(\Omega\cap U)$, so that $w$ extends continuously to the closure of $\Omega\cap U$. Assume that there are positive constants $c_5$, $c_7$, $C'$, $M_4$ and $M_5$ such that
\[
 c_7d\le w_0\le C'd,\qquad-C'\le\mathcal{L}w_0\le-c_5\qquad\text{and}\qquad|w|\le M_5\qquad\text{in }D,
\]
that $\mathcal{L}w\le M_4$ in $\Omega\cap U$ in the sense of distributions, and that $w\ge0$ on $\Delta_{2\rho}(P')$. Then there are $\rho_*>0$ and $c_6>0$, depending only on $M$, $\kappa$, $c_7$ and $C'$, such that for $0<\rho\le\rho_*$
\[
 w\ge-C_{10}\,w_0\qquad\text{in }\Omega_{\rho_1/c}(P'),\qquad\text{where}\qquad C_{10}=\frac{cM_5}{c_6\rho}+\frac{M_4}{c_5}.
\]
\end{lemma}
\begin{proof}
Let $\tau=s-s_*$, where $s_*$ is the smallest value of $s$ on $D$, so that $\mathcal{L}\tau=-1$ and $0\le\tau\le8\rho^2$. Let $Z_0$ and $Z_1$ be the solutions of $\mathcal{L}Z=0$ in $D$ given by \cref{lem:dirichlet}\,(i) that agree on the Kolmogorov boundary with $w_0$ and $w^-$, and let $Z_2=\tau-Z_3$, where $Z_3$ is the solution with boundary values $\tau$. Then $\mathcal{L}Z_2=-1$, the function $Z_2$ vanishes on the Kolmogorov boundary, and $0\le Z_2\le\tau$ by \cref{lem:dirichlet}. Since $w_0/c_5-Z_2$ is a supersolution that is nonnegative on the Kolmogorov boundary, $w_0$ is a supersolution and $w_0-C'\tau$ is a subsolution, \cref{lem:dirichlet} gives
\[
 c_5Z_2\le w_0\qquad\text{and}\qquad Z_0\le w_0\le Z_0+8C'\rho^2\qquad\text{in }D.
\]
Moreover $Z_0>0$ in $D$ by the strong maximum principle of Bony~\cite{Bon69}, because $Z_0\ge0$, the face $\{x=x_{P'}+8M\rho\}$ of $D$ is reached from every point of $D$ along $\partial_x$, and $Z_0=w_0>0$ on this face, because $w_0\ge c_7d$ and $d\ge4M\rho$ there by \eqref{eq:LNlip}.

We claim that $w\ge-Z_1-M_4Z_2$ in $D$. Let $\eta>0$ be small, and let $w^\varepsilon$ be a mollification of $w$ at a scale $\varepsilon$ so small that the $2\varepsilon$-neighborhood of the closure of $D\cap\{d>\eta\}$ lies in $\Omega\cap U$. Such $\varepsilon$ exist, because $\overline D$ is a compact subset of $U$ and $d$ is continuous. Since $w$ is Lipschitz continuous in $\Omega\cap U$ with bounded $w_{xx}$, and since $\mathcal{L}w\le M_4$ there in the sense of distributions, the function $w^\varepsilon$ is smooth on a neighborhood of the closure of $D\cap\{d>\eta\}$, the commutators of $\mathcal{L}$ with the mollification are bounded by $M_6\varepsilon$ there, where $M_6$ may depend on $w$, and $|w^\varepsilon-w|\le M_6\varepsilon$. Since $\mathcal{L}\tau=-1$, it follows that
\[
 \mathcal{L}\bigl(w^\varepsilon+Z_1+M_4Z_2+M_6\varepsilon\tau\bigr)\le0\qquad\text{in }D\cap\{d>\eta\}.
\]
On the part of the Kolmogorov boundary of $D\cap\{d>\eta\}$ where $d=\eta$, the function in parentheses is at least $w^\varepsilon\ge-M_6(\varepsilon+\eta)$ after enlarging $M_6$, because $w\ge0$ on the graph. The remaining part lies in the Kolmogorov boundary of $D$, where $Z_1=w^-$ and $Z_2=0$, so that the function is at least $w^+-M_6\varepsilon$ there. \Cref{lem:dirichlet}\,(ii), applied in $D\cap\{d>\eta\}$, gives the lower bound $-M_6(\varepsilon+\eta)$ there. Letting $\varepsilon\searrow0$ and then $\eta\searrow0$ proves the claim.

By \eqref{eq:refpoints}, the reference points $\Xi^\pm_{\rho_1}(P')$ satisfy $d\ge\rho_1/C$ with $C=C(M)$. Hence
\[
 Z_0\bigl(\Xi^\pm_{\rho_1}(P')\bigr)\ge w_0\bigl(\Xi^\pm_{\rho_1}(P')\bigr)-8C'\rho^2\ge\frac{c_7\rho}{Cc_0^2}-8C'\rho^2\ge c_6\rho\qquad\text{with}\qquad c_6=\frac{c_7}{2Cc_0^2}
\]
for $\rho\le\rho_*=c_7/(16CC'c_0^2)$. Since $Z_1\le M_5$, applying \cref{lem:LN} to $Z_1+\varepsilon Z_0$ and $Z_0$ and letting $\varepsilon\searrow0$, we obtain
\[
 Z_1\le\frac{cM_5}{c_6\rho}\,Z_0\le\frac{cM_5}{c_6\rho}\,w_0\qquad\text{in }\Omega_{\rho_1/c}(P').
\]
Together with $M_4Z_2\le(M_4/c_5)\,w_0$ and the claim, this proves the lemma.
\end{proof}

\begin{proposition}\label{prop:asian-cone}
For every compact set $E\subset\{(t,A):\ 0<t<T,\ A>K\}$, the cone estimate \eqref{eq:cone} holds for the Asian problem near the part of the free boundary that lies over $E$.
\end{proposition}
\begin{proof}
Here \eqref{eq:cone} reads $|u_I|+|u_t|\le C_0u_S$. Let $E'$ be a compact neighborhood of $E$ in $\{0<t<T,\ A>K\}$. We use \cref{lem:asian-boundary} for $E'$, together with the set $N$ and the constants $L_1$, $M$ and $M_2$ from its proof, and we write $d=S-\Phi(I,t)$. By \cref{lem:asian-shift}, $-u_I\le L_1u_S$ and $u_t\le L_1u_S$ in $N\cap\{u>0\}$. It remains to bound $u_t$ from below at the points of $\{u>0\}$ near the free boundary over $E$.

At these points we have the following three properties. First, $u$ is convex in $S$ and vanishes at $S=\Phi(I,t)$, so \cref{lem:asian-boundary} gives $u_S\ge u/d\ge c_3d$ for $d\le\eta_0$, and $u_S\le Md$ because $u_S$ vanishes on $\Gamma$. Second, since $[S\partial_S,\mathcal{K}]=S\partial_I$ and $\mathcal{K}u=f$ in $\{u>0\}$, we have
\[
 \mathcal{K}(Su_S)=Sf_S-Su_I=-S\Bigl(u_I+\frac1t\Bigr),
\]
and this lies between $-S/t$ and $-S/(2t)$ when $d$ is small, because $0\le-u_I\le L_1Md$. Third, for $h>0$ the function
\[
 v_h(t,S,I)=\frac{u(t+h,S,I)-u(t,S,I)}{h}
\]
satisfies $|v_h|\le M_2$, and $v_h\ge0$ on $\Gamma$. The coefficients of $\mathcal{K}$ do not depend on $t$, and $\mathcal{K}u\le f$ almost everywhere by \cref{lem:asian-incentive}. Hence $\mathcal{K}v_h\le M_3$ almost everywhere in $\{u>0\}$ near $\Gamma$, where $M_3$ bounds $|f_t|$ on $N$. We keep $h$ fixed until the end of the proof. Then $v_h$ is Lipschitz continuous and $(v_h)_{SS}$ is bounded, so this inequality also holds in the sense of distributions in $\{u>0\}$ near $\Gamma$.

As in the proof of \cref{prop:C1gamma}, we reverse time, change $a$ away from the region under consideration, and put $w_0=e^{-\varphi}Su_S$ and $w_h=e^{-\varphi}v_h$. Then $\mathcal{L}w_0=e^{-\varphi}\mathcal{K}(Su_S)-cw_0$ and $\mathcal{L}w_h=e^{-\varphi}\mathcal{K}v_h-cw_h$. By the first property, $w_0$ is comparable to $d$, and since $w_h$ is bounded, the second and third properties give $-C'\le\mathcal{L}w_0\le-c_5$ and $\mathcal{L}w_h\le M_4$ near $\Gamma$, the latter in the sense of distributions. Moreover $|w_h|\le M_5$ with $M_5=M_2\max e^{-\varphi}$, and $w_h\ge0$ on $\Gamma$. None of these constants depends on $h$. By \cref{lem:asian-boundary} for $E'$, the bound \eqref{eq:halfLip} holds near the image of every point of the free boundary over $E$, so the sets $\Omega_\rho(P')$ of \eqref{eq:LNsets} are available there, and finitely many of these local constructions cover the free boundary over $E$. Fix a graph point $P'$ near the image of such a point, and let $\rho\le\rho_*$ be so small that the closure of $\Omega_{2\rho}(P')$ lies in an open neighborhood $U$ of the free boundary in which the bounds above hold on $\{u>0\}\cap U$. There $w_h$ is Lipschitz continuous with bounded $(w_h)_{xx}$, because $v_h$ has these properties, so the hypotheses of \cref{lem:single-scale} are satisfied. \Cref{lem:single-scale}, applied to $w_0$ and $w=w_h$, gives $w_h\ge-C_{10}w_0$ in $\Omega_{\rho_1/c}(P')$ with $\rho_1=\rho/c_0^2$, where $C_{10}$ does not depend on $h$ and on $P'$. By \eqref{eq:LNballs} and \eqref{eq:LNsets}, this set contains the points $P'\circ(x,0,0)$ with $0<x<\rho_1/c^2$. Undoing the conjugation and the time reversal and letting $h\searrow0$, we obtain $u_t\ge-C_{10}Su_S$ at the points of $\{u>0\}$ near the free boundary over $E$ where $d<\rho_1/c^2$. At the other points near this part of the free boundary, the first property and the convexity of $u$ in $S$ give $u_S\ge c_3\min\{\rho_1/c^2,\eta_0\}$, while $|u_t|\le M_2$. This proves the proposition.
\end{proof}

\begin{proof}[Proof of \cref{cor:asian}]
By \cref{lem:asian-incentive,lem:asian-monotone,lem:asian-gap}, for $0<t<T$ and $A>K$ the point $(t,S,tA)$ lies in the exercise region if and only if $S\le S_*(t,A)$, where $S_*(t,A)=\Phi(tA,t)>0$. Fix such a point $(t,A)$ and a compact neighborhood $E$ of it in $\{0<t<T,\ A>K\}$. By \cref{lem:asian-incentive,lem:asian-monotone,lem:asian-gap,prop:asian-cone}, there is an open set $Q$ containing the free boundary over the interior of $E$ in which $u$ is a solution of \eqref{eq:obstacle} under \eqref{eq:H} that satisfies \eqref{eq:cone}. By \cref{thm:A}\,(iii), the function $\Phi$ is of class $C^{1,\gamma}$ near the free boundary over the interior of $E$ for some $\gamma\in(0,1)$. Since $S_*(t,A)=\Phi(tA,t)$, the same holds for $S_*$.
\end{proof}

The constants of this section degenerate as $A\searrow K$, as $t\searrow0$ and as $t\nearrow T$, and \cref{cor:asian} says nothing about the behavior of $S_*$ there.
\begin{remark}\label{rem:asian-limits}
As $A\searrow K$ the critical price tends to zero. Indeed, stopping at $T$ gives
\[
 V(t,S,tA)\ge e^{-r(T-t)}\,\mathbb{E}\Bigl[\Bigl(\frac{SJ_T}{T}-\frac{K(T-t)}{T}\Bigr)^+\Bigr]=:g(t,S)\qquad\text{for }A\ge K,
\]
with $J_T$ as in the proof of \cref{lem:asian-monotone}, and $g(t,S)$ is positive and nondecreasing in $S$ and tends to zero as $S\searrow0$. Exercise at $(t,S,tA)$ requires $V=A-K$, so that $g(t,S_*(t,A))\le A-K$, and hence $S_*(t,A)\to0$ as $A\searrow K$. At maturity one expects $S_*(t,A)\to A+rT(A-K)$ as $t\nearrow T$. Indeed, by the computation in the proof of \cref{lem:asian-incentive}, the discounted payoff $e^{-rt}(A_t-K)$ has the drift $-e^{-rt}f$, so that immediate exercise is favorable near maturity where $f>0$, and $\{f>0\}=\{S<A+rt(A-K)\}$. We have not proved this limit, and the limit $t\searrow0$ involves the initial state. We do not pursue these limits. Lookback options and floating-strike Asian options lead to obstacle problems of the same type, and the comparisons of \cref{lem:asian-shift} have analogues for them, but we have not checked the details.
\end{remark}

\medskip
\noindent\textbf{Funding.} Junkee Jeon is supported by the National Research Foundation of Korea (NRF) grant funded by the Korea government [Grant No. RS-2026-25591924]. Takwon Kim is supported by the National Research Foundation of Korea (NRF) grant funded by the Korea government [Grant No. RS-2024-00351151]. Jehan Oh is supported by the National Research Foundation of Korea (NRF) grant funded by the Korea government [Grant Nos. RS-2025-00555316 and RS-2025-25415411].

\noindent\textbf{Use of artificial intelligence.} Large language model tools were used for language editing and for checking routine computations. The authors verified every statement and proof. The authors assume responsibility for all content.

\end{document}